\RequirePackage{tikz}
\documentclass[sn-mathphys]{sn-jnl}

\usepackage{subcaption}

\jyear{2021}%

\theoremstyle{thmstyleone}%
\theoremstyle{thmstyletwo}%

\theoremstyle{thmstylethree}%

\begin{document}

\title[Modelling of reversible tissue electroporation]{Modelling of reversible tissue electroporation and its thermal effects in drug delivery}


\author*[1]{\fnm{Nilay} \sur{Mondal}}\email{nilay.mondal@iitg.ac.in}
\equalcont{These authors contributed equally to this work.}

\author[1]{\fnm{D. C.} \sur{Dalal}}\email{durga@iitg.ac.in}
\equalcont{These authors contributed equally to this work.}


\affil*[1]{\orgdiv{{Department of Mathematics}, \orgname{Indian Institute of Technology Guwahati}, \orgaddress{\street{Amingaon}, \city{Guwahati}, \postcode{781039}, \state{Assam}, \country{India}}}}

%


\abstract{Electroporation is a very useful tool for drug delivery into various diseased tissues of the human body. This technique helps to improve the clinical treatment by transferring drugs into the targeted cells rapidly. In electroporation, drug particles enter easily into the intracellular compartment through the temporarily permeabilized cell membrane due to the applied electric field. In this work, a  mathematical model of drug delivery focusing on reversible tissue electroporation is presented. In addition, the thermal effects on the tissue, which is an outcome of Joule heating, are also considered. This model introduces a time-dependent mass transfer coefficient, which is significant to drug transport.  Multiple pulses with low voltage are applied to reach sufficient drugs into the targeted cells.  According to the physical circumstances, a set of differential equations are considered and solved. The changes in drug concentration with different parameters (e.g., diffusion coefficient, drug permeability,  pulse length, and pulse number) are analyzed. The model optimizes the electroporation parameters to uptake drugs sufficiently into the cells with no thermal damage. This model can be used in clinical experiments to predict drug uptake into the infected cells by controlling the model parameters according to the nature of infections.}

\keywords{Electroporation, Drug Delivery, Reversibly Electroporated Cells, Low Voltage Multiple Pulses, Temperature}



\maketitle

\section{Introduction}
Now-a-days, various diseases like tumors, cancer, etc., are treated using gene therapy and chemotherapy. Different physical approaches are developed to deliver drugs or genes into the diseased tissues. New techniques like electroporation, micro-injection, laser, ultrasound are frequently used in physical delivery systems \cite{Bolhassani2011}. Electroporation is a biophysical phenomenon in which cell membranes are temporarily destabilized due to the application of external electric fields, and during this destabilizing period, nano-meter size pores are formed in lipid bilayer of the cell membrane \cite{Krassowska2007,Miklavcic2014,Miklavcic2016,Miklavcic2020}.  During electroporation, the transmembrane potential ($V_m$) increases due to the application of electric field and cell membrane gets destabilized \cite{Pavlin2003,Kotnik2019}. Neumann et al. \cite{Neumann1989} presented a mathematical model to evaluate the transmembrane potential of a spherical cell in a uniform electric field. Different molecules such as drugs, ions, dyes, antibodies, oligonucleotides, RNA and DNA, etc., are transported into the intracellular space through the temporary pores \cite{Satkauskas}.
Electroporation has vast applications in both \textit{in vitro} and \textit{in vivo} \cite{Miklavcic2018Feb,Kotnik2019,Miklavcic2020}. Some tissues such as Skin, liver, tumors, and muscles are generally targeted through \textit{in vivo} electroporation \cite{Widera2000}.

Experimental and theoretical studies indicate that two types of electroporation occur, namely reversible electroporation and irreversible electroporation. In reversible electroporation, cell membranes are temporarily permeabilized, and the treated cells survive the process. If the applied electric field is very high, electroporation leads to cell death due to permanent permeabilization of the cell membrane, and the process is called irreversible electroporation \cite{Rubinsky2007,Miklavcic2012,Miklavcic2018SSept}. Reversibly electroporated cells go through the treatment, and cell membranes start the resealing process after completion of the pulse period \cite{Davalos2003,Granot}. Irreversibly electroporated cells do not survive the treatment, and pores do not reseal \cite{Jiang2015}. It is seen that the formation of pores is rapid (in microseconds); on the other hand, the membrane resealing process takes a very long time (in minutes). During the late 1970s, the cell membrane's resealing effect was found experimentally \cite{Kinosita,Gauger} and in 1997, Saulis \cite{Saulis1997} investigated it theoretically.

The drug transport into the targeted cells depends on the physical properties of the tissue like conductivity, size, shape, and distribution of cells, along with the electrical parameters like number of pulses, pulse strength, and pulse duration \cite{Pavlin,Pucihar2011}. In the various experiment, short and long duration pulses are generally used for cell electroporation. An inhomogeneous short duration electric pulse is applied in single-cell electroporation. 
It is seen that the pore density in the cell membrane increases with the application of a short duration (1$\mu$s - 100 $\mu$s) high voltage electric field  \cite{Weaver2003,Pavlin}. To electroporate a tissue, bulk electroporation is generally used by applying homogeneous long duration (100 $\mu$s - 100 ms) pulses \cite{Granot,Becker2013,Miklavcic2018Feb}. Joule heating and electrophoresis are the two critical outcomes due to the application of long duration pulse. In Joule heating,  localized temperature increases, and it helps to raise permeability of cell membrane through lipid chain melting  \cite{Smith,Becker2016}. It is necessary to emphasize on the increment of permeability and control of cell death for the drug uptake in reversibly electroporated cells. Efforts are being made continuously to improve mathematical models so that the above concerns can be addressed. A brief overview of the most important theoretical models on drug delivery with tissue electroporation is presented below.

In the area of electroporation based drug delivery, Granot and Rubinsky \cite{Granot}  proposed a mass transfer model in 2008 to deliver drug into the cells with reversible electroporation. In their study, the pore creation model \cite{Krassowska2007} is used and how the rate of mass transfer to the cells increases with pore density is described.
In 2014,  Kalamiza et al. \cite{Kalamiza20141950}  established a dual-porosity model focusing on solute diffusion in the biological tissue after electroporation. This is probably the first study that built a relationship between the increment of resultant permeability and the macroscopic drug transport into the cells. The mass transfer coefficient (MTC) that depends on the increased membrane permeability due to the formation of electropores was proposed in the model. The effects of irreversible electroporation is not taken into account. Later in 2016, Boyd and Becker \cite{Bradley2016} used the previous model to drug delivery using tissue electroporation. In their study, mass transport equations for extracellular and intracellular spaces are provided separately. Pore resealing effects on drug uptake is also studied. Argus et al. \cite{Becker2017} modified the last model in 2017 by considering irreversible electroporation into account. The study represents a three-equation mass transport model to calculate drug concentrations in extracellular space, reversibly and irreversibly electroporated cells. Recently in the study \cite{Mondal2021}, a mathematical model for drug delivery is presented in which both reversible and irreversible electroporations with the application of uniform electric field are considered. Two separate approaches: drug delivery during and after electroporation are described in the model. Only temporal changes of drug concentrations are emphasized but spatial changes of drug concentrations are not considered in this study.

In the present study, a mathematical model of drug delivery is developed considering bulk electroporation on a tissue. Only reversible electroporation is emphasized in order to cure whole infected tissue by delivering drug as a medicine. On this account, low voltage multiple pulses are used to avoid cell death and permanent permeabilization of the membrane. The principle goal of this study is to analyze the drug transport phenomena in extracellular and intracellular domains of the targeted tissue, when drug is administrated from one side. The novelty of this study is the investigation of the effects of thermal changes in the tissue on cell survival.  Due to the application of repeated pulses in electroporation, a time dependent MTC is considered. The MTC is basically a function of pore density and it plays a significant role to control drug transport from extracellular to intracellular space. Some differential equations along with appropriate initial and boundary conditions are prescribed to represent the governing biological circumstances. 
Analytical and numerical methods are used to solve the differential equations wherever applicable. In numerical simulation,  equations are discretized by finite difference method (FDM). 
The parameters' dependencies on mass transfer are studied in details. It is possible to introduce a required amount of drugs into the targeted cells of a tissue through this model. No cell damage or cell death is noticed due to thermal effect as low voltage pulses are applied. Qualitative comparisons between the present study and the existing theoretical and experimental studies \cite{Wong2011,Kuldeep2021,Kalamiza20141950} are made in order to authenticate the advocated model.
The numerical results of this model can be helpful in the pharmaceutical and medical research.

\section{Problem Formulation}

\begin{figure}[h!]
	\centering
	\begin{tikzpicture}[scale=0.4]
	
	\draw[fill,blue!30] (-6,-6) rectangle (8,8);
	\draw[ultra thick] (-6,-6) rectangle (8,8);
	
	\draw[fill,gray!70](0,0) circle (1cm);
	\draw[ultra thick, dashed] (0,0) circle (1cm);
	
	\draw[fill,gray!70] (-2,2) circle (1cm);
	\draw[ultra thick,dashed] (-2,2) circle (1cm);
	
	\draw[fill,gray!70](3,3) circle (1cm);
	\draw[ultra thick, dashed] (3,3) circle (1cm);
	
	\draw[fill,gray!70] (-4,4) circle (1cm);
	\draw[ultra thick,dashed] (-4,4) circle (1cm);

	\draw[fill,gray!70](4,-3) circle (1cm);
	\draw[ultra thick, dashed] (4,-3) circle (1cm);
	
	\draw[fill,gray!70] (4,6) circle (1cm);
	\draw[ultra thick,dashed] (4,6) circle (1cm);
	
	\draw[fill,gray!70](-4,-4) circle (1cm);
	\draw[ultra thick, dashed] (-4,-4) circle (1cm);
	
	\draw[fill,gray!70] (-2,-2) circle (1cm);
	\draw[ultra thick,dashed] (-2,-2) circle (1cm);
	
	\draw[fill,gray!70] (-4.5,0) circle (1cm);
	\draw[ultra thick,dashed] (-4.5,0) circle (1cm);
	
	\draw[fill,gray!70] (-2,6) circle (1cm);
	\draw[ultra thick,dashed] (-2,6) circle (1cm);
	
	\draw[fill,gray!70] (5,1) circle (1cm);
	\draw[ultra thick,dashed] (5,1) circle (1cm);
	
	\draw[fill,gray!70] (0,4) circle (1cm);
	\draw[ultra thick,dashed] (0,4) circle (1cm);
	
	\draw[fill,gray!70] (0,-4) circle (1cm);
	\draw[ultra thick,dashed] (0,-4) circle (1cm);
	
	\draw[thick, <-] (6,-5) -- (10,-5);
	\draw (13.5,-5) node[]{Extracellular Space};
	\draw (5,-5) node[]{$C_E$};
	
	\draw[thick, <-] (5.5,0.5) -- (10,0.5);
	\draw (12,0.5) node{Reversibly };
	\draw (12,-0.3) node{electroporated cell};
	\draw (5,1.2) node{$C_{RE}$};
	
	\draw[thick, <-] (5,6) -- (10,6);
	\draw (13,6) node{Cell membrane };
	
	\draw[ultra thick, ->] (-10,6) -- (-6,6);
	\draw[ultra thick, ->] (-10,1) -- (-6,1);
	\draw[ultra thick, ->] (-10,-4) -- (-6,-4);
	\draw (-8,0) node{Drug};
	
	\end{tikzpicture}
	
	\caption{A schematic diagram of injecting drug  into a biological tissue.}
	\label{2tissue1}	
\end{figure}
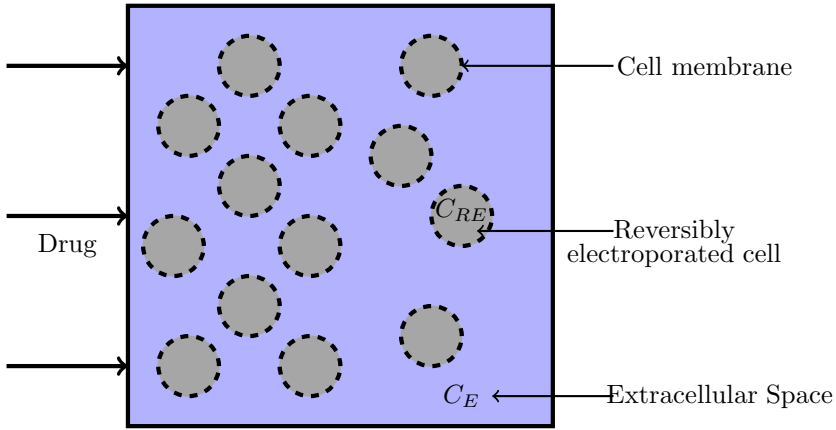
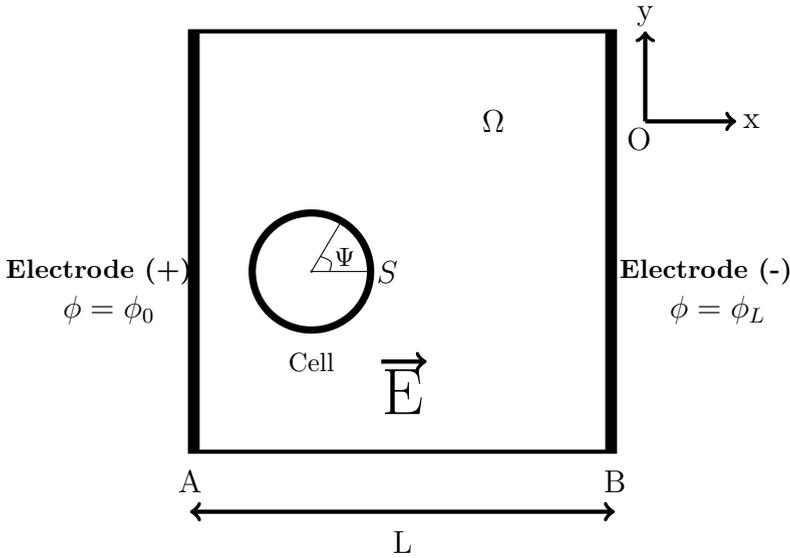
\begin{figure}[h!]
	\centering
	\begin{tikzpicture}[scale=0.4]
	
	\draw[fill,black!100] (-6,-6) rectangle (-5.7,8);
	\draw[fill,black!100] (7.7,-6) rectangle (8,8);
	\draw[ultra thick] (-6,-6) rectangle (8,8);
	\draw[ultra thick] (-2,0) circle (1.9cm);
	\draw[ultra thick] (-2,0) circle (2.0cm);
	\draw[fill,black] (-2,0) circle (0.02cm);
	\draw[] (-2,0) -- (0,0);
	\draw[] (-2,0) -- (-1,1.7);
	\draw[] (-1.4,0) .. controls (-1.3,0.35) .. (-1.7,0.5);
	\draw (-0.9,0.5) node[]{$\Psi$};
	\draw (-2,-3) node[]{Cell};
	
	\draw (-6,-7) node[]{\large A};
	\draw (8,-7) node[]{\large B};
	
	\draw[ultra thick,<->](-6,-8) -- (8,-8);
	\draw (1,-9) node[]{\large L};
	\draw[ultra thick,->](0.3,-3) -- (1.8,-3);
	\draw (1,-4) node[]{\Huge E};

	\draw (11,0) node[]{\bf Electrode (-)};
	\draw (9.5,-0.5) node[anchor=north west]{\large $\phi=\phi_L$};
	
	\draw (-9,0) node[]{\bf Electrode (+)};
	\draw (-10.5,-0.5) node[anchor=north west]{\large $\phi=\phi_0$};
	
	\draw[ultra thick,->](9,5) -- (9,8);
	\draw[ultra thick,->](9,5) -- (12,5);
	\draw (12.5,5) node[]{\large x};
	\draw (9,8.5) node[]{\large y};
	\draw (8.8,4.4) node[]{\large O};
	
	\draw (4,5) node[]{\large $\Omega$};
	\draw (0.5,0) node[]{\large $S$};
	
	\end{tikzpicture}
	
	\caption{A schematic representation of bulk electroporation on a square tissue using two electrodes at A and B.}
	\label{2electrode1}
\end{figure}

This study considers a tissue, which is assumed to be cubical having side length  $L$. The tissue domain consists of extracellular and intracellular spaces. Both areas are separated by the porous cell membrane, which controls the mass transport from the extracellular to the intracellular domain. A schematic diagram (shown by Fig. \ref{2tissue1}) is provided to visualize the complete tissue structure as well as the drug injecting process into the tissue. In order to electroporate the whole tissue, two electrodes with different potential values $\phi_0$ and $\phi_L$ are placed along two different boundaries, respectively (as shown in Fig. \ref{2electrode1}). A uniform electric field $E$ is induced in the whole tissue directed from the positive to the negative electrode. Here, bulk electroporation is applied to electroporate all the cells of the tissue using multiple pulses with low voltage. The transmembrane potential ($V_m$) is increased, and nanometer-sized pores are formed in the cell membrane due to application of electric pulses. It is considered that the strength of the applied electric field does not exceed the value 28 V mm$^{-1}$ to maintain the critical transmembrane potential ($\approxeq$ 1 V). The maximum pores are created at the poles $\Psi=0$, $\pi$ as transmembrane potentials are maximum at those locations. No new pores are formed at the poles $\Psi=\frac{\pi}{2}$, $\frac{3\pi}{2}$ as the transmembrane potential is zero at those poles.  All cells in the tissue are electroporated reversibly due to application of low voltage pulses. After electroporation, the cell membrane gets permeabilized and allows drugs to enter the cell. The resealing process in the electroporated cells happens during the rest time between the pulses, and the drug transport occurs only in this rest time. The drug uptake into cells is affected due to the membrane resealing effect.

\section{Model Development and Mathematical Equations}
Based on the arrangement of two parallel electrodes, a uniform electric field is generated inside the tissue. The potential ($\phi$) distribution inside the tissue is obtained by solving the  Laplace equation \cite{Sel2005,Corovic,Bradley2016}
\begin{eqnarray}\label{2eq1}
\nabla^2 \phi=0,
\end{eqnarray}
with boundary conditions: $\phi(x=0)=\phi_0$ and $\phi(x=L)=\phi_L$.\\
The uniform electric field ($E$) in the domain $\Omega$ is calculated by the magnitude of the gradient of the potential and mathematically  expressed as, 
\begin{eqnarray}\label{2E}
E= \lvert \vec{\nabla} \phi \lvert
\end{eqnarray}
The transmembrane potential  $V_m $ due to the effect of electric field $E$ is obtained by the relation \cite{Krassowska1999}
\begin{equation}\label{transm_potential}
V_m=1.5E\times r_c\cos \Psi ,
\end{equation}
where $r_c$ is the radius of the cell and  $\Psi$ is the angle between the direction of electric field  and the normal to the cell membrane at the position where $V_m$ is calculated.

The pore density $N$ in the cell membrane is obtained from the time dependent equation proposed by Krassowka and Filev \cite{Krassowska2007} as,
\begin{equation}\label{2pore_density}
\frac{dN}{dt}=\alpha A \left[1- \frac{N}{N_0}A^{-q}\right],
\end{equation}
where $A=\exp \left[\left( \frac{V_m}{V_{ep}} \right)^2\right]$, t is the time, $\alpha$ is the pore creation rate coefficient,  $V_{ep}$ the characteristic voltage of electroporation, $N_0$ the equilibrium pore density at $V_m = 0$ and $q$ is an electroporation constant. Initially, it is considered that there is no pore in the membrane i.e., $N(0)=0$. 

The total  pore area $(A_P)$  of a single cell is calculated as, $A_P=\pi R_P^2\cdot N_P$, where $R_P$ is the pore radius and $N_P$ is the total number of pore in a cell membrane.
As cells are reversibly electroporated, the cell membrane goes to reseal after electroporation and the pore area starts decreasing with time. The modified pore area can be mathematically expressed as an exponentially decreasing function of time \cite{Granot},
\begin{equation}\label{2eq5}
A_P=\pi R_P^2\cdot N_P(t_{ep})\cdot \exp\left(-\frac{t}{\tau}\right),
\end{equation}
where $t_{ep}$ is the pulse duration, and $\tau$ is the resealing time constant.	

\subsection{Drug transport phenomenon in the tissue}
In this model, the drug  is introduced as a medicine at the left boundary of the targeted tissue. It is assumed that the drug is injected before the electroporation process starts. Multiple pulses are applied for a certain time duration ($t_{ep}$) mentioning a long time gap ($t_M$) between the pulses. Only in this time gap, mass transfer occurs from extracellular to intracellular compartment. The electric pulses are applied periodically until a sufficient drug reaches into the cells. The mass transfer rate depends on the total pore area of the cell. The pore area is a function of time due to the membrane resealing effect. So, the MTC ($\mu$) is a function of time and it is mathematically defined as,
\begin{equation}\label{2mut}
\mu(t)=\left(\frac{4\pi^2 R^{2}_Pr_c^2 }{V_0}\right)P\cdot N(t_{ep})\cdot \exp\left(-\frac{t}{\tau}\right),
\end{equation}
where $V_0$ is the volume of a cube containing a cell and $P$ is the permeability of drug particles  across the  cell membrane.

The drug diffuses from left to right through the extracellular space of the tissue, and when pores are formed, drug enters into the cells.  The drug concentrations in the extracellular space and in the electroporated cells are obtained by the mass transport equations and mathematically can be expressed as,
\begin{eqnarray}\label{2eq7}
\frac{\partial C_E}{\partial t}=\vec{\nabla}.\left(D \vec{\nabla} C_E\right) - \left(\frac{1-\varepsilon}{\varepsilon}\right)\mu (t)\times
\left(C_E - C_{RE}\right),
\end{eqnarray} 
\begin{equation}\label{2eq8}
\frac{\partial C_{RE}}{\partial t}=\mu (t)\times (C_E - C_{RE}),
\end{equation}
where  $C_E$ and $C_{RE}$ are the drug concentrations in the extracellular space and reversibly electroporated cells respectively; $D$ is the effective diffusion coefficient of the drug in the extracellular space; $\varepsilon$ is the porosity and it is defined as the volumetric ratio of extracellular volume to the total volume of the tissue. The above equations are solved considering appropriate initial and boundary conditions, which are as follows:\\
\begin{equation}\label{2eq8.1}
\left.\begin{split}
&C_E(x, y, 0)=C_1(x,y),\qquad C_{RE}(x, y, 0)=0,\\
&C_E(0, y, t)=C_2,\qquad \qquad
\left(\frac{\partial C_E}{\partial x}\right)_{x=L}=0,\\
&\left(\frac{\partial C_E}{\partial y}\right)_{y=0}=0,\qquad \qquad
\left(\frac{\partial C_E}{\partial y}\right)_{y=L}=0,
\end{split}\qquad\right\}
\end{equation}
where 
\begin{equation*}
C_1(x,y)=\left \{\begin{split}
&C_2;\qquad x=0,~~ 0\leq y\leq L,\\
&0;\qquad \qquad\text{otherwise}
\end{split}\qquad\right.
\end{equation*}
\subsection{Thermal effects in the tissue}
\begin{figure}[h!]
	\centering
	\begin{tikzpicture}[scale=0.45]
	\draw[fill,red!20] (-9.5,-9.5) rectangle (11.5,11.5);
	\draw (1,10) node[]{\large Outside tissue};
	\draw (1,8.6) node[]{\large $T_b$};
	\draw (1,-7.4) node[]{\large $T_b$};
	\draw (-7,1) node[]{\large $T_b$};
	\draw (9,1) node[]{\large $T_b$};
	
	\draw[fill,blue!20] (-6,-6) rectangle (8,8);
	\draw[ultra thick] (-6,-6) rectangle (8,8);
	\draw (1,-2) node[]{\large Tissue};
	\draw (1,1) node[]{\large $T$};
	
	\draw (-6,-7) node[]{\large $x=0$};
	\draw (8,-7) node[]{\large $x=L$};
	
	\draw (-7.8,-6) node[]{\large $y=0$};
	\draw (-7.8,8) node[]{\large $y=L$};
	
	
	\draw[ultra thick,->](14,5) -- (14,8);
	\draw[ultra thick,->](14,5) -- (17,5);
	\draw (15.2,4.3) node[]{\large x};
	\draw (13.4,6) node[]{\large y};
	
	\end{tikzpicture}
	
	\caption{A schematic representation of temperature distribution in a tissue due to electroporation.}
	\label{2temp1}
\end{figure}
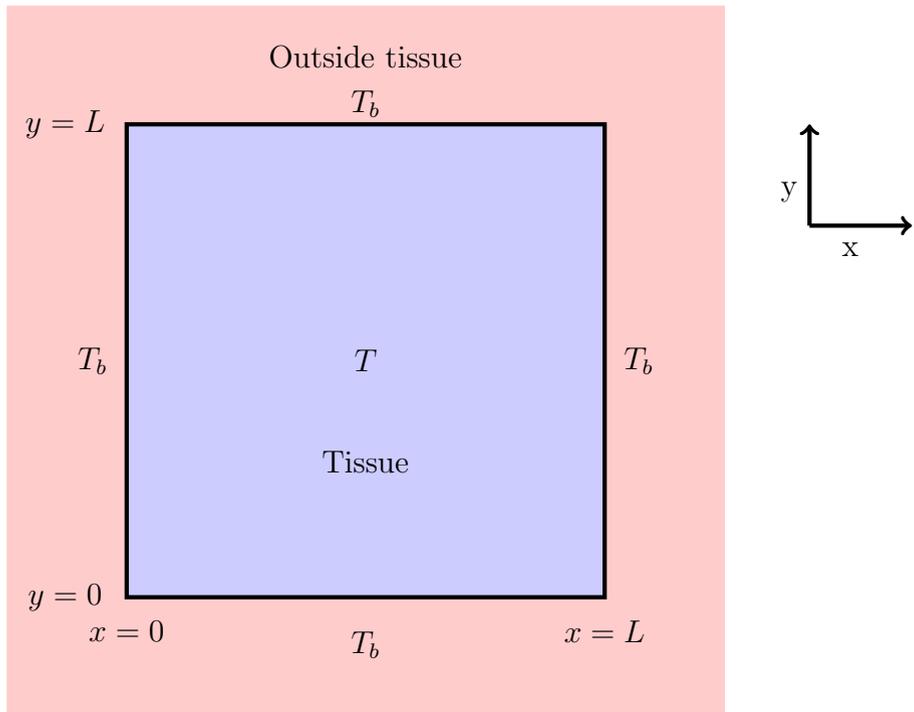
The membrane structure of a cell has the electrical resistivity even if the lipid layers of the membrane are thin. The cell membranes act as capacitors with high resistance \cite{Foster1989}. Moreover, when the tissue gets exposed to an external electric field, the cell membranes become charged up. Elevated current density due to tissue electroporation guides to Joule heating, which is significant especially when longer electric pulses are applied. The Joule heating may not be biologically significant when electric pulses are applied for short duration, but in order to investigate the underlying mechanisms behind electro-chemotherapy, the role played by the elevated temperature is important. Furthermore, the study of temperature variation is essential because it helps in formation of stable permeable structures and hence plays a crucial role in electrically mediated molecular transport \cite{Pliquett2003}.  Also, the increment in temperature in cells plays a significant role in the dielectric breakdown of cell membrane. It is seen from experiment that the potential difference for dielectric breakdown decreases with the increase in temperature \cite{Coster1975a}.  The mathematical formulation of Joule effect is represented through the following equations.

The change in temperature due to ongoing electroporation inside the tissue is obtained by the heat equation \cite{Garcia2014},
\begin{equation}\label{2eqtemp1}
\rho  c  \frac{\partial T}{\partial t}=k \nabla^2 T + Q_J,
\end{equation}
with initial condition:  $T(x, y, 0)=T_b$,\\
and boundary conditions:\\
\begin{center}
	$\left(\frac{\partial T}{\partial x}\right)_{x=0}=\frac{h}{k}(T-T_b),~~~$
	$\left(\frac{\partial T}{\partial x}\right)_{x=L}=-\frac{h}{k}(T-T_b)$,\\
	$\left(\frac{\partial T}{\partial y}\right)_{y=0}=\frac{h}{k}(T-T_b),~~~$
	$\left(\frac{\partial T}{\partial y}\right)_{y=L}=-\frac{h}{k}(T-T_b)$,\\
\end{center}
where $\rho$ is the density, $c$, the specific heat and $k$ is the thermal conductivity of the tissue; $T_b ~(\leq T$), the body temperature; $h$, the heat transfer coefficient.\\
The resistive heating term $Q_J$ refers to the heat source due to Joule heating as a result of applied electric pulse and it can be defined as
\begin{align}
Q_J=\sigma \lvert \vec{\nabla} \phi \lvert^2,
\end{align}
where $\sigma$ is the tissue electrical conductivity.

The temperature is expected to increase continuously with time over the time duration of electroporation and on the termination electric pulse, temperature may decrease with time. The governing temperature equation for this process may be the standard heat equation
\begin{equation}\label{2eqtemp2}
\rho  c  \frac{\partial T}{\partial t}=k \nabla^2 T,
\end{equation}
with initial condition:  $T(x, y, 0)=T_{ep}$,\\
and boundary conditions:
\begin{center}
	$\left(\frac{\partial T}{\partial x}\right)_{x=0}=\frac{h}{k}(T-T_b),~~~$
	$\left(\frac{\partial T}{\partial x}\right)_{x=L}=-\frac{h}{k}(T-T_b)$,\\
	$\left(\frac{\partial T}{\partial y}\right)_{y=0}=\frac{h}{k}(T-T_b),~~~$
	$\left(\frac{\partial T}{\partial y}\right)_{y=L}=-\frac{h}{k}(T-T_b)$.\\
\end{center}
Here, $T_{ep}$ is the temperature of the tissue after application of electric pulse over the time duration $t_{ep}$.
\section{Method of Solutions}
\begin{table}[h!]
	\caption{The parameter's detail used for simulation of the model:}
	\vspace{0.2cm}
	\label{2tab1}
	\centering
	\begin{tabular}{llll}
		\hline
		Symbol & Value  & Definition & Source   \\
		\hline 
		&&&\\
		$\sigma$ & 0.241 S m$^{-1}$ & Tissue electrical conductivity & \cite{Becker2017}  \\
		$r_c$ & 25 $\mu$m & Cell radius & \cite{Kalamiza20141950}  \\
		$\alpha$  & $10^9$ m$^{-2}$ s$^{-1}$ & Pore creation coefficient  & \cite{Krassowska2007}  \\
		$V_{ep}$ & 0.258 V& Characteristic voltage & \cite{Krassowska2007}  \\
		$N_0$ & $1.5\times10^9$ m$^{-2}$ & Equilibrium pore density   & \cite{Krassowska2007} \\
		$q$ & $2.46$ & Electroporation constant & \cite{Granot} \\
		$D$ & $10^{-4}$ mm$^2$ s$^{-1}$ & Effective diffusion coefficient &\\
		$R_P$& 0.8 nm& Pore radius & \cite{Granot} \\
		$\varepsilon$& 0.18 & Porosity & \cite{Becker2017} \\
		$P$ & $5\times 10^{-4}$ mm s$^{-1}$ & Permeability of drug & \cite{Granot} \\
		$E$ & 28 V mm$^{-1}$ & Electrical field & \\
		$C_2$ & 1 M & Initial drug concentration &\\
		$L$ & 1 mm  & Length of the square & Fig. \ref{2electrode1} \\
		$\phi_0$ & 28 V  & Potential along $x=0$ & Fig. \ref{2electrode1}\\
		$\phi_L$ & 0 V  & Potential along $x=L$ & Fig. \ref{2electrode1}\\
		$\rho$&1060 kg m$^{-3}$&Density of the tissue & \cite{Garcia2014}\\
		$c$&3600 J kg$^{-1}$ K$^{-1}$ & Specific heat capacity of the tissue & \cite{Garcia2014}\\
		$k$&0.502 W m$^{-1}$ K$^{-1}$ & Thermal conductivity of the tissue & \cite{Garcia2014}\\
		$h$&50 W m$^{-2}$ K$^{-1}$ & Heat transfer coefficient of the tissue & \cite{Davalos2003} \\
		$T_b$&37 $^\circ$ C = 310.15 K & Normal body temperature & \cite{Garcia2014}\\
		$t_{ep}$ & 80 ms & Pulse length (ON TIME) & \\
		$t_{M}$ & 600 s & Time for mass transfer (OFF TIME) & \\
		$PN$ & 10 & Pulse number &\\
		&&&\\
		\hline
	\end{tabular}
\end{table}

The governing equations are solved  analytically wherever possible otherwise, numerical methods are applied to obtain approximate solutions. Since, two electrodes are placed vertically along left and right boundaries with different potentials, a uniform electrical field is induced over the whole tissue region. Due to parallel electrode arrangement (see Fig. \ref{2electrode1}), it may be considered that $\phi$ is independent of $y$ i.e. $\frac{\partial \phi}{\partial y}=0$.   Eq. \eqref{2eq1} along with the boundary conditions is solved analytically in a square region $\Omega$ ($0\leq x \leq L$, $0\leq y \leq L$) and the potential distribution is obtained as,
\begin{equation}
\phi(x)=\frac{(\phi_L-\phi_0)}{L}x+\phi_0.
\end{equation}
Hence, from the Eq. (\ref{2E}), a uniform electric field in the tissue region $\Omega$ can be obtained as, 
\begin{equation*}
E=\frac{(\phi_0-\phi_L)}{L}.
\end{equation*}
Eq. \eqref{2pore_density} is also solved analytically considering $\Psi = 0$ or $\pi$ and the initial condition $N(0)=0$.  The pore density $N(t)$ is found as,
\begin{equation}
N(t)=N_0A^q\left[1- \exp \left(-\frac{\alpha t}{N_0 \cdot A^{q-1}}\right)   \right].
\end{equation}
The Eq. \eqref{2mut} gives the MTC as,
\begin{equation}\label{2soln_mut}
\mu(t)=\left(\frac{4\pi^2 R^{2}_Pr_c^2 P N_0A^q}{V_0}\right)\cdot \left[1- \exp \left(-\frac{\alpha t}{N_0 \cdot A^{q-1}}\right)   \right].
\end{equation}

The Eqs. \eqref{2eq7}-\eqref{2eqtemp2} are solved numerically using  finite difference method. 
In numerical computation, the Eqs. \eqref{2eq7}-\eqref{2eq8} are discretized using the forward time centered space (FTCS) scheme. The discretized equations are as follows,
\begin{align}\label{2eq15}
\left(C_E\right)_{i,j}^{n+1}=&a\left(C_E\right)_{i+1,j}^{n}+a\left(C_E\right)_{i-1,j}^{n}+b(t)\left(C_E\right)_{i,j}^{n} \\ \nonumber
&+c\left(C_E\right)_{i,j+1}^{n}+c\left(C_E\right)_{i,j-1}^{n}+d(t)\left(C_{RE}\right)_{i,j}^{n}, \quad 
\begin{array}{c}
i=1,2,\cdots, M_1\\
j=1,2,\cdots, M_2
\end{array}
\end{align}
\begin{align}\label{2eq16}
\left(C_{RE}\right)_{i,j}^{n+1}=\left(C_{RE}\right)_{i,j}^{n}+\mu(t)\Delta t\left[\left(C_E\right)_{i,j}^{n}-\left(C_{RE}\right)_{i,j}^{n}\right], \quad \begin{array}{c}
i=1,2,\cdots, M_1 \\
j=1,2,\cdots, M_2
\end{array}
\end{align}
where $a=\frac{D\Delta t}{(\Delta x)^2}$, 
$b(t)=1-\left[2D\left(\frac{1}{(\Delta x)^2}+\frac{1}{(\Delta y)^2}\right)+\frac{1-\varepsilon}{\varepsilon}\mu(t)\right]\Delta t$, 
$c=\frac{D\Delta t}{(\Delta y)^2}$, and $d(t)=\frac{1-\varepsilon}{\varepsilon}\mu(t)\Delta t$. Here, $\Delta t$ is the time step size, $\Delta x$ and $\Delta y$ are the step sizes for space, $M_1$ and $M_2$ are the numbers of grid along $x$-axis and $y$-axis respectively.\\
The stability condition used to solve the Eqs. \eqref{2eq15}-\eqref{2eq16} is
\begin{align*}
\Delta t<\frac{1}{2}\times \frac{(\Delta x)^2 (\Delta y)^2}{D[(\Delta x)^2+ (\Delta y)^2]}.
\end{align*}
Here, $M_1=101$, $M_2=101$, $\Delta x=0.01$, $\Delta y=0.01$ and $\Delta t=0.2$ are chosen for simulation.

The Eqs. \eqref{2eqtemp1}-\eqref{2eqtemp2} are also discretized by FTCS scheme and the discretized equations are as follows,
\begin{align}\label{2eq17}
T_{i,j}^{n+1}=eT_{i+1,j}^{n}+eT_{i-1,j}^{n}+fT_{i,j}^{n}
+gT_{i,j+1}^{n}+gT_{i,j-1}^{n}+lQ_J, \quad 
\begin{array}{c}
i=1,2,\cdots, M_1 \\
j=1,2,\cdots, M_2
\end{array}
\end{align}
\begin{align}\label{2eq18}
T_{i,j}^{n+1}=eT_{i+1,j}^{n}+eT_{i-1,j}^{n}+fT_{i,j}^{n}
+gT_{i,j+1}^{n}+gT_{i,j-1}^{n} \quad 
\begin{array}{c}
i=1,2,\cdots, M_1 \\
j=1,2,\cdots, M_2
\end{array}
\end{align}
where $e=\frac{k\Delta t}{\rho c (\Delta x)^2}$, $f=1-\frac{2k}{\rho c}\left[\frac{1}{(\Delta x)^2}+\frac{1}{(\Delta y)^2}\right]\Delta t$, $g=\frac{k\Delta t}{\rho c (\Delta y)^2}$, and $l=\frac{\sigma\Delta t}{\rho c}$.\\
The stability condition used for the Eqs. \eqref{2eq17}-\eqref{2eq18} is
\begin{align*}
\Delta t<\frac{1}{2}\times \frac{(\Delta x)^2 (\Delta y)^2}{\frac{k}{\rho c}[(\Delta x)^2+ (\Delta y)^2]}.
\end{align*}
Here, $M_1=101$, $M_2=101$, $\Delta x=0.01$, $\Delta y=0.01$ and $\Delta t=0.00002$ are chosen for simulation. Two different $\Delta t$ are used to solve the  Eqs. \eqref{2eq15}-\eqref{2eq16} and  \eqref{2eq17}-\eqref{2eq18} separately. The details of the parameters used in the model are shown in  Table \ref{2tab1} and the suitable values of the parameters are taken from the existing literature. 

\section{Results and Discussion}
\begin{figure}[h!]
	\centering
	\begin{tikzpicture}[scale=0.35]
	
	\draw[ultra thick] (-6,-6) rectangle (8,8);
	
	\draw[ultra thick, ->] (-10,6) -- (-6,6);
	\draw[ultra thick, ->] (-10,2) -- (-6,2);
	\draw[ultra thick, ->] (-10,-2) -- (-6,-2);
	\draw (-8,4) node{\large Drug};
	\draw (-8.5,0) node{$C_E=1$ M};
	\draw (1,2) node{$C_E=0$ M};
	\draw (1,0) node{$C_{RE}=0$ M};
	\draw (1,9) node{$\frac{\partial C_E}{\partial y}=0$ };
	\draw (1,-7) node{$\frac{\partial C_E}{\partial y}=0$ };
	\draw (10,0) node{$\frac{\partial C_E}{\partial x}=0$ };
	\end{tikzpicture}
	
	\caption{A schematic diagram of the model equations, which are presented in the Eqs. \eqref{2eq7}-\eqref{2eq8.1} for injecting drug  into a biological tissue.}
	\label{2model}	
\end{figure}
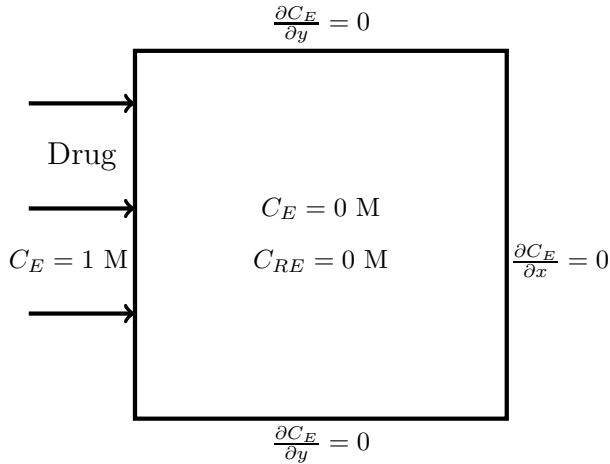
\begin{figure}[h!]
	\centering
	\begin{tabular}{cc}
		\includegraphics[width=.5\linewidth]{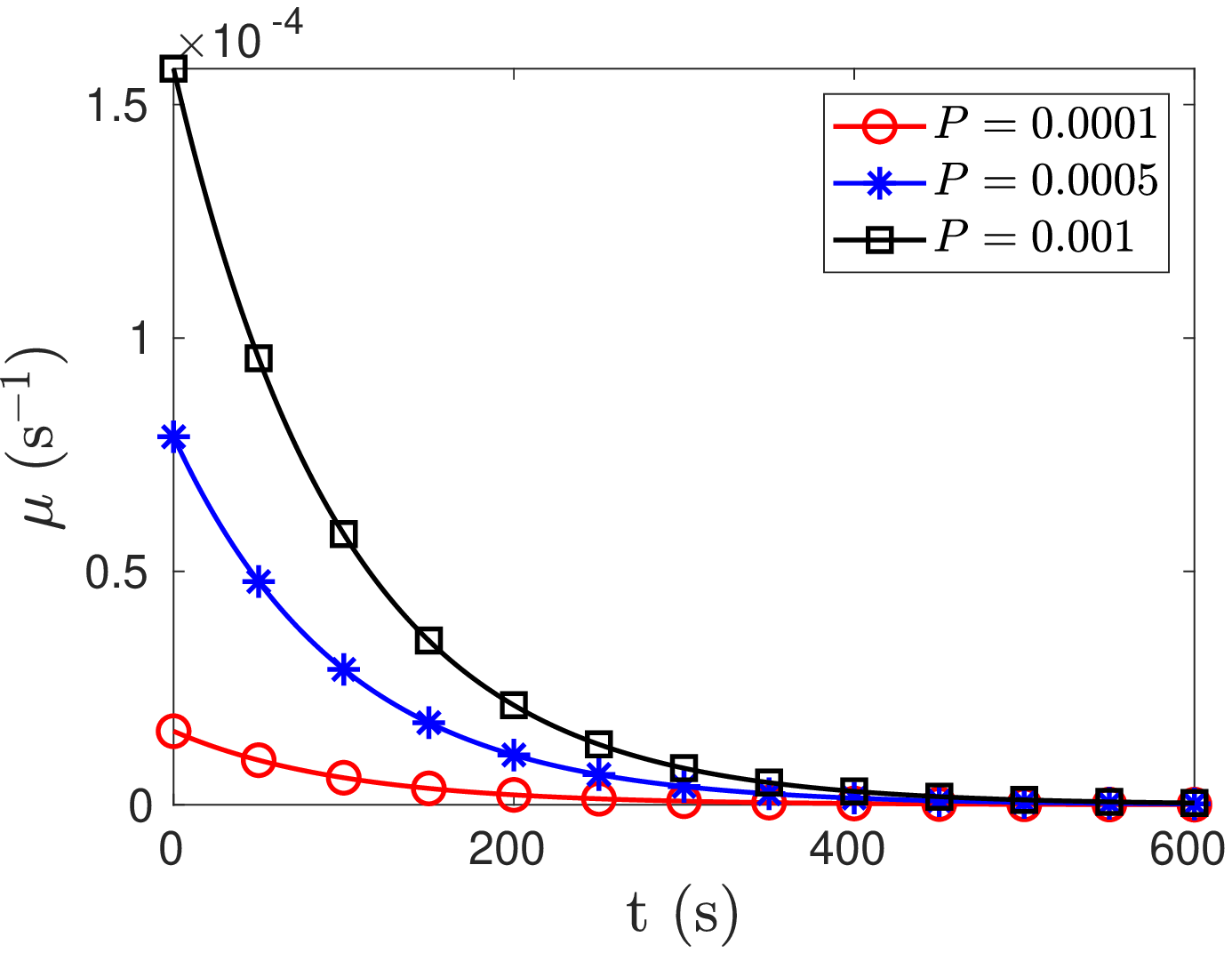}
		&
		\includegraphics[width=.5\linewidth]{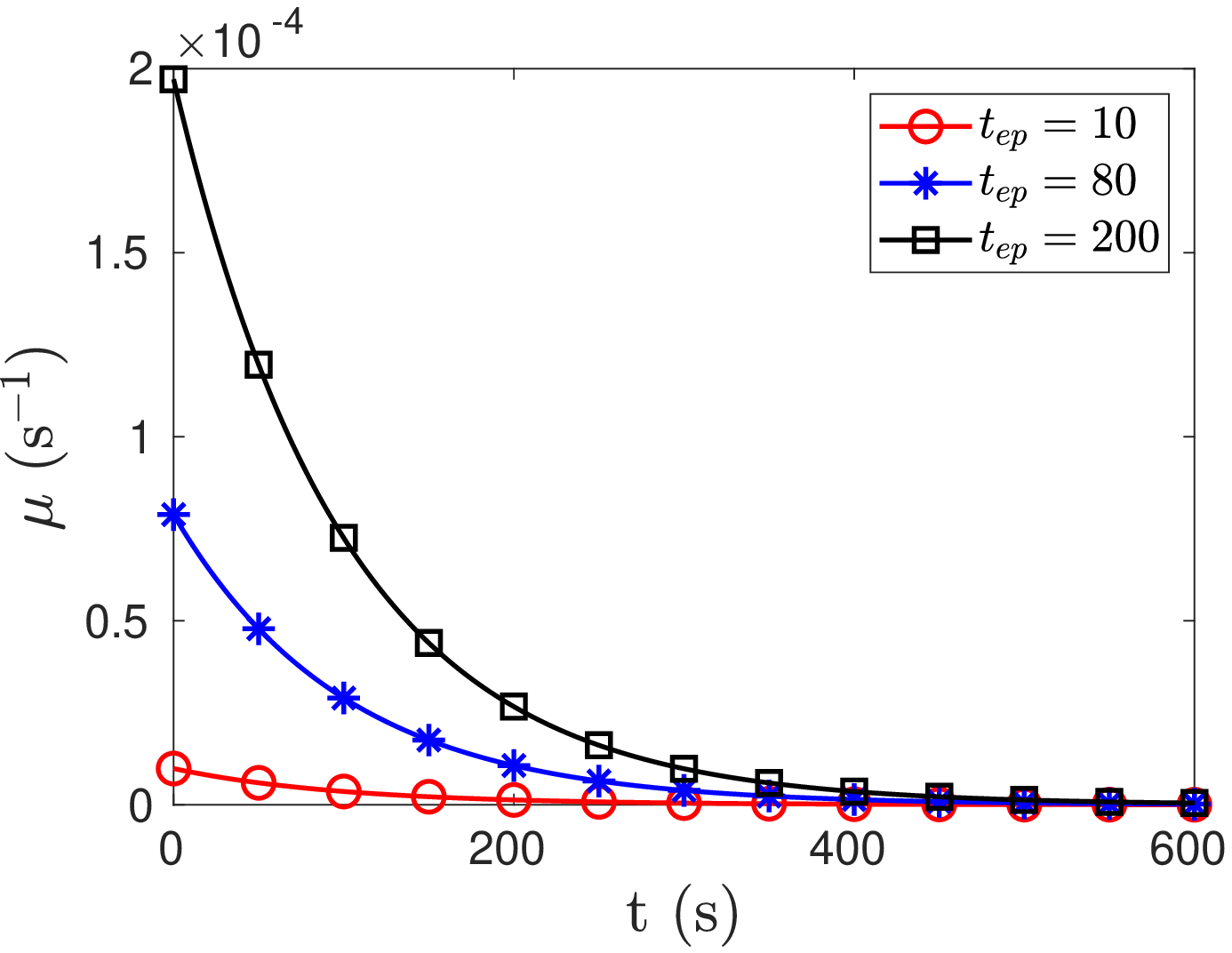}
		\\(a) & (b) \\	
	\end{tabular}
	\caption{The changes of MTC with time after electroporation for different (a) permeability $P$ (mm s$^{-1}$) and (b) pulse length $t_{ep}$ (ms).}
	\label{2mu}
\end{figure}

In the present study, the drug transport phenomenon in the diseased tissue is represented by a mathematical model. The main aim of this study is to focus on drug penetration through a targeted tissue when drug is injected from one side. The physical overview of the drug transport in this model is shown in as a diagram in Fig. \ref{2model}. A qualitative analysis is performed to visualize the drug dynamics into the electroporated tissue through graphical representations in Figs. \ref{2mu} - \ref{2temp_surface}. The detailed discussion on the plotted graphs of drug concentrations for different values of drug permeability ($P$), diffusion coefficient ($D$) and pulse length ($t_{ep}$) is made.

Fig. \ref{2mu} describes how the MTC ($\mu$) changes with time after  electroporation for different permeability and various pulse lengths. In Fig. \ref{2mu}a, it is noticed that  $\mu$ increases with drug permeability. In Fig. \ref{2mu}b, it is also seen that the parameter $\mu$ increases with the pulse length because of the higher pore density for the application of longer pulses. The physical observation may be explained through the mathematical expression given by the Eq. \eqref{2mut}.  In both the graphs of Fig. \ref{2mu}, the MTC decreases with time due to membrane resealing. The increased MTC helps to increase the drug uptake rate into the cells from the extracellular region.

\subsection{Drug transport phenomenon in the tissue}
This section presents a full overview of the drug transport process in the targeted tissue when the drug is applied from the left side. The process includes continuous molecular diffusion in the homogeneous extracellular fluid medium, and drug uptake into the intracellular domain from the extracellular domain through the  permeabilized cell membrane.

\begin{figure}[h!]
	\centering
	\begin{tabular}{cc}
		\includegraphics[width=.5\linewidth]{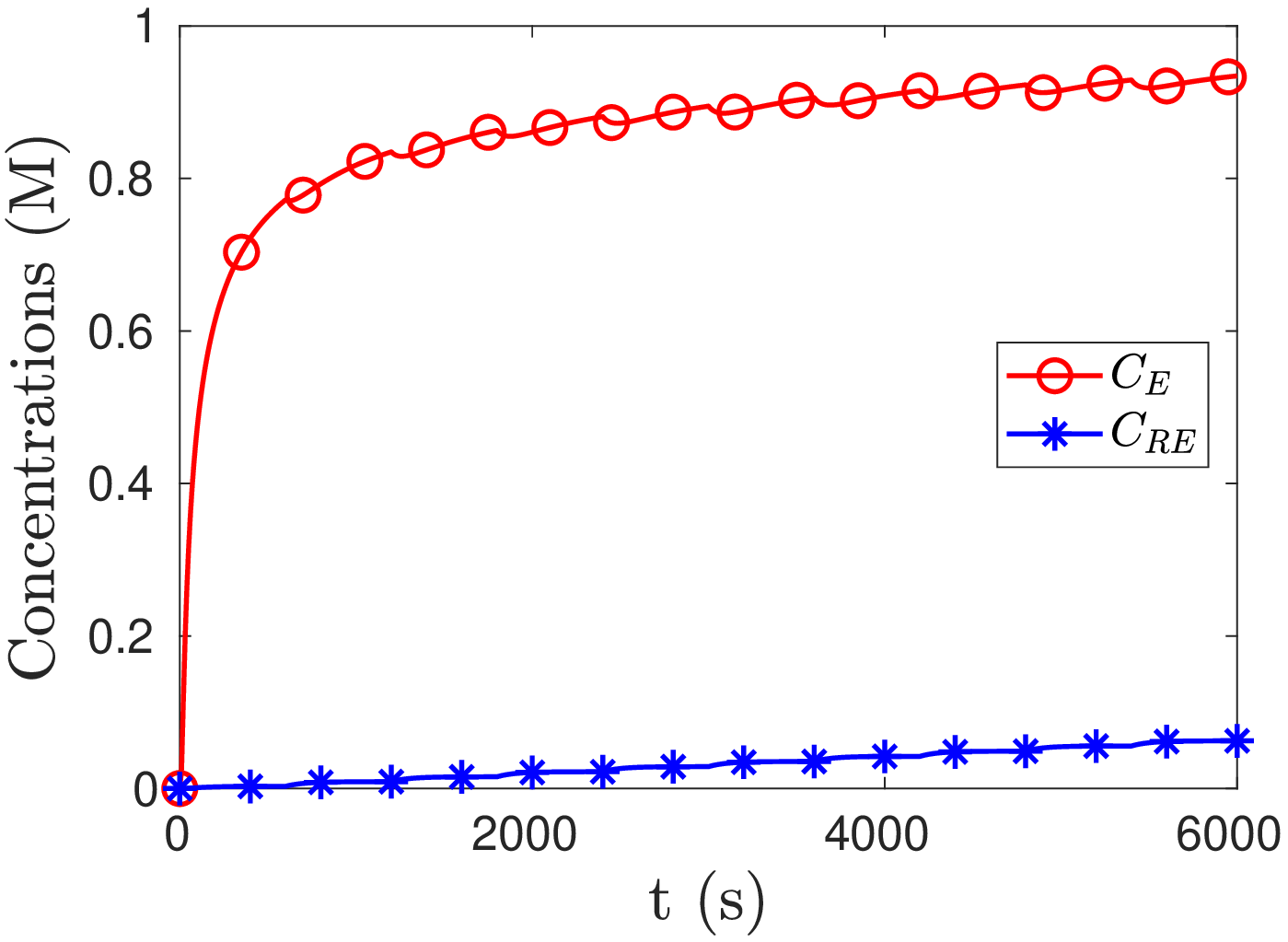}
		&
		\includegraphics[width=.5\linewidth]{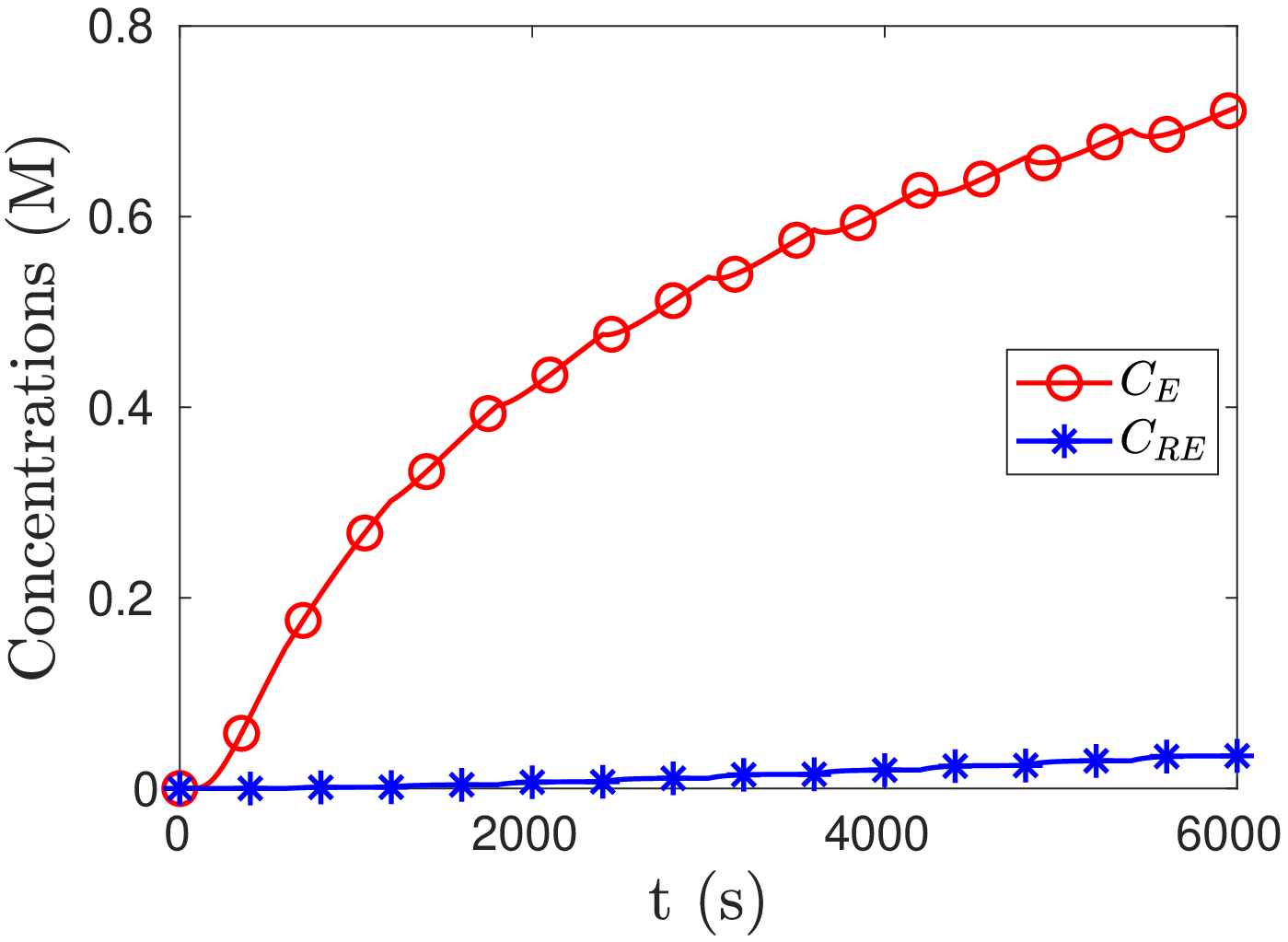}
		\\(a) $x = 0.1$ mm & (b) $x = 0.5$ mm
	\end{tabular}
	\begin{tabular}{c}
		\includegraphics[width=.5\linewidth]{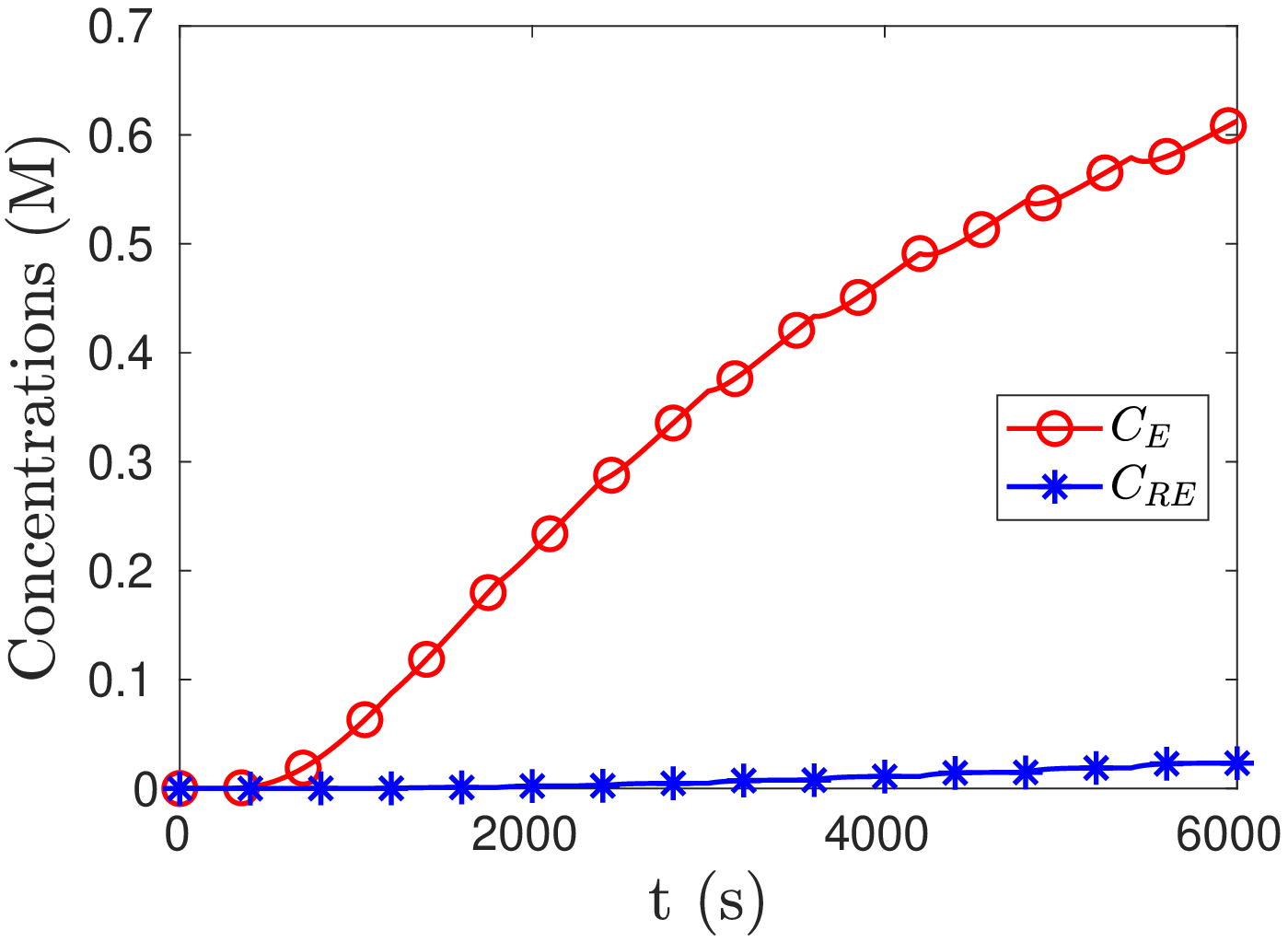}
		\\(c) $x = 0.9$ mm	
	\end{tabular}
	\caption{The drug concentration ($C_E$, $C_{RE}$) profiles, which are changing with time, at three different locations parallel to x-axis along the line $y = 0.5$.  Here, $C_E$, $C_{RE}$ are calculated for $P = 0.0005$ mm s$^{-1}$ and $D = 0.0001$ mm$^2$ s$^{-1}$.}
	\label{2cons}
\end{figure}
Fig. \ref{2cons} shows the time dependent drug concentrations in extra- and intra-cellular spaces at the locations $x = 0.1$, $x = 0.5$ and $x=0.9$. From Fig. \ref{2cons}a, it is observed that the drug concentration at $x = 0.1$ in the extracellular region increases within a short time. This is because the location is near to the drug source at the left boundary so the molecular diffusion to that point occurs quickly. On the other hand, the drug concentration in the intracellular region increases slowly at that time due to the low permeability of the cell membrane as fewer pores are formed immediately after the electroporation. It is also noticed that both the profiles overlap after some time, near about 150  s. Because, the membrane permeability increases due to sufficient number of pores formed by this time (near 150  s from the start of electroporation) and the dug uptake into the cells occurs faster. In Fig. \ref{2cons}b, it is shown that the drug concentrations ($C_E$, $C_{RE}$) at $x = 0.5$ for the entire time are relatively less than the same at $x = 0.1$. As the point ($0.5, 0.5$) is in the middle and far from the drug source, less amount of drug reaches there. Another reason is that the drug uptake in the cells near the boundary is higher. Similar observation is noticed through the Fig. \ref{2cons}c. As the point (0.9, 0.5) is far from the drug source, relatively less amount of drug reaches in the intracellular as well as extracellular region.
\subsubsection{Comparison}
\begin{figure}[h!]
	\centering
	\includegraphics[width=7cm]{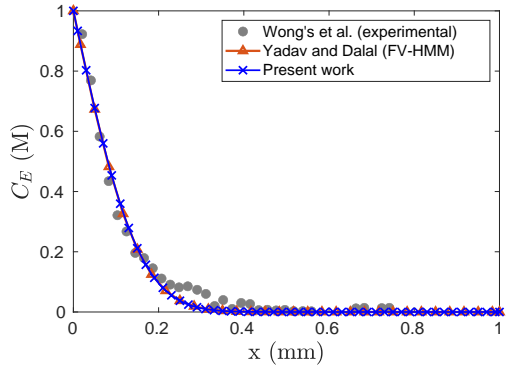}	
	\caption{The present study is compared with the experimental work and the FV-HMM model.}
	\label{2comparison1}
\end{figure}
\begin{figure}[h!]
	\centering
	\begin{tabular}{cc}
		\includegraphics[width=.5\linewidth]{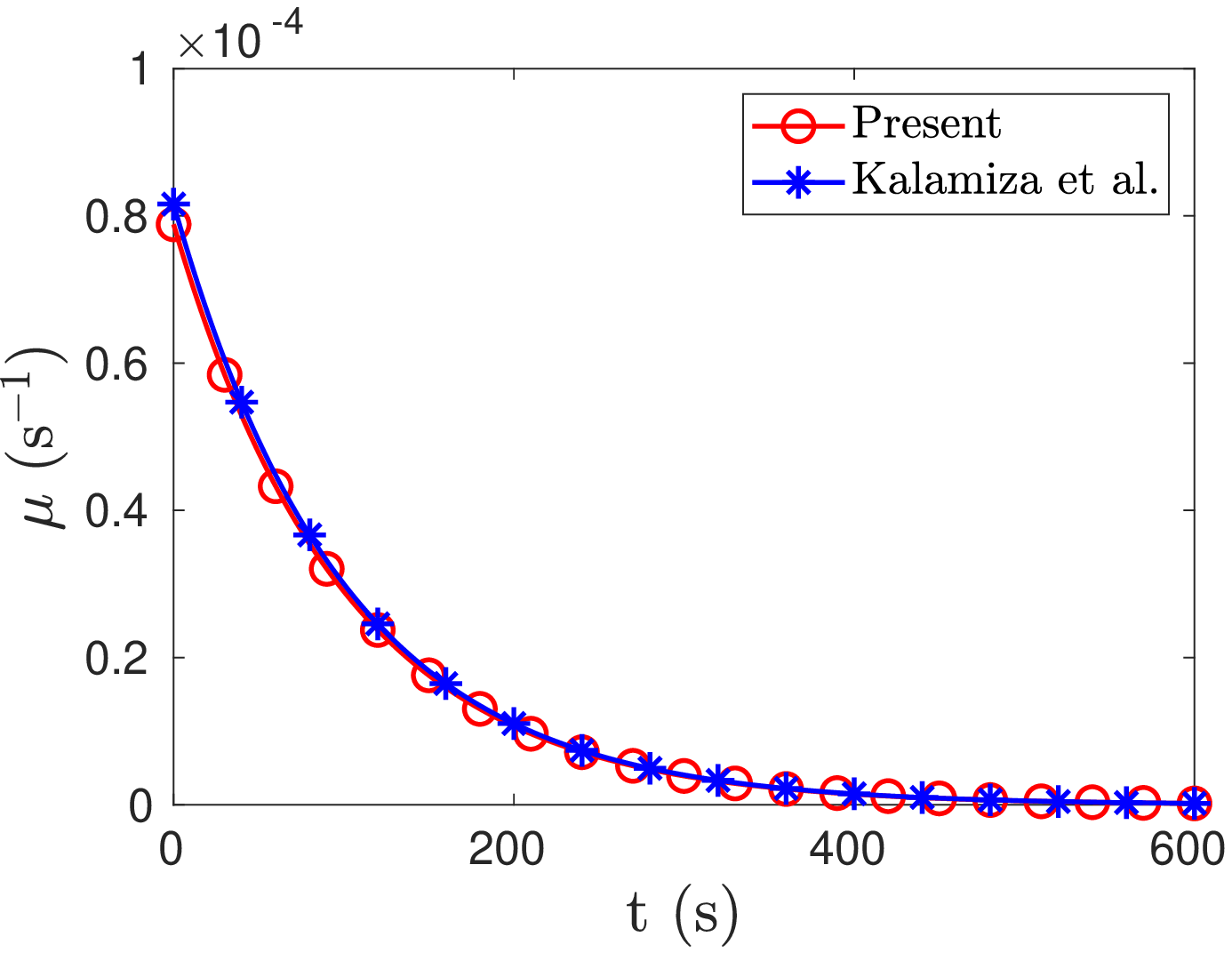}
		&
		\includegraphics[width=.5\linewidth]{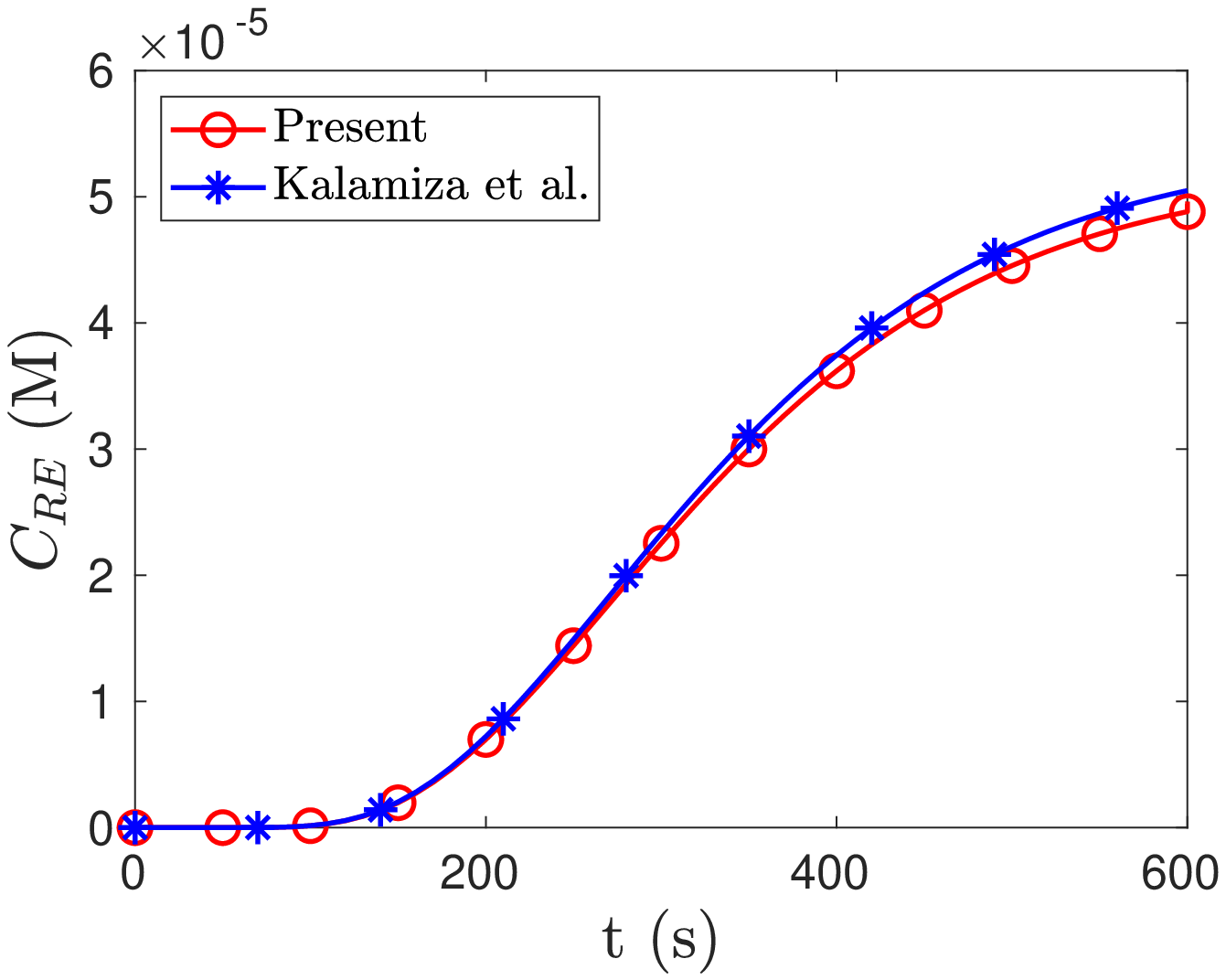}
		\\(a) & (b) \\		
	\end{tabular}
	\caption{Comparison between the present results with the existing model. Here, the concentration ($C_{RE}$) are calculated at (0.5, 0.5) for $P=0.0005$ mm s$^{-1}$, $D=0.0001$ mm$^2$ s$^{-1}$ and $t_{ep}=80$ ms.}
	\label{2comparison2}
\end{figure}
To validate the present model, two comparisons are made. In first comparison, electroporation is not taken into account and results are compared with the experimental results of Wong's et al \cite{Wong2011} and the theoretical results of Yadav and Dalal \cite{Kuldeep2021}.  As  no electroporation is considered, $\mu=0$ is taken to solve the concentration equations in the present model and the comparison is shown in the Fig. \ref{2comparison1}. 
In the experiment \cite{Wong2011}, a test on nanoparticle penetration through a collagen gel (compared to tissue media) is performed. In FV-HMM model, the drug penetration through a tissue is studied and the heterogeneous multiscale method is used to solve the model. In all the three models, the drug penetration in the tissue is studied for a duration of 12 hours.  The Fig. \ref{2comparison1} shows that the drug penetration (present study) in the extracellular media along x-axis is very close to those in the other studies. 

Another comparison is made when electroporation is incorporated i.e. $\mu\ne0$. In this case, the present results are compared with the work of Kalamiza et al. \cite{Kalamiza20141950}. The electroporation based MTC is taken from the model \cite{Kalamiza20141950}. In order to compare the results, the MTC that includes the pore resealing effect can be written as \cite{Kalamiza20141950},
\begin{equation}\label{2mut_kalamiza}
\mu(t)=\frac{3Df_p}{d_mr_c}\cdot \exp\left(-\frac{t}{\tau}\right)=8.16\times 10^{-5}\cdot \exp\left(-\frac{t}{\tau}\right),
\end{equation}
where $f_p=3.4\times 10^{-8}$ and $d_m=5\times 10^{-6}$ mm are considered. 

The results for both MTCs (as in Eqs. \eqref{2mut} and \eqref{2mut_kalamiza}) are plotted in  the Fig. \ref{2comparison2}.  From the plots, it can be observed that there is a good agreement between the results.
\subsubsection{Effects of drug permeability on drug concentrations}
\begin{figure}[h!]
	\centering
	\begin{tabular}{cc}
		\includegraphics[width=.5\linewidth]{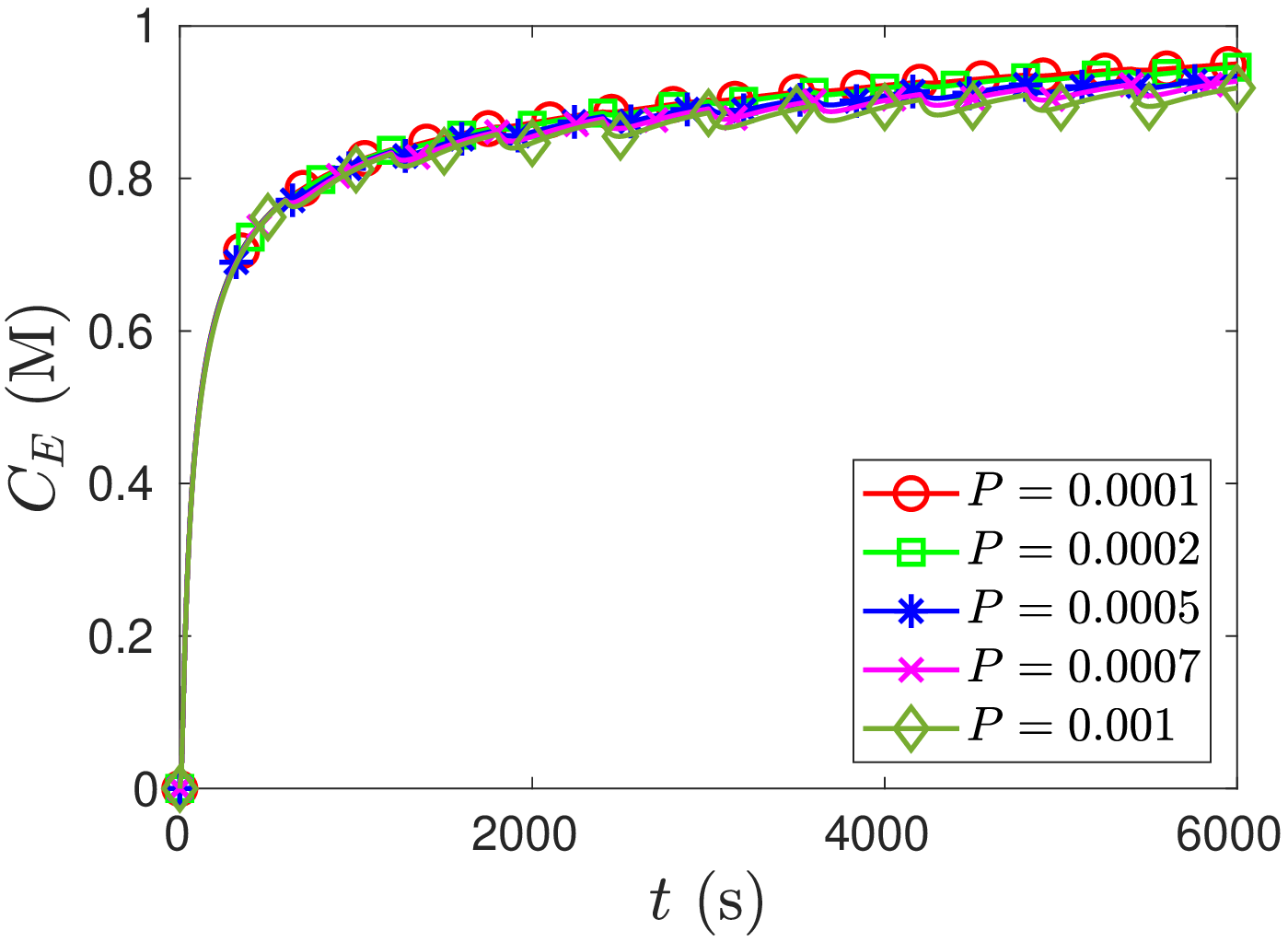}
		&
		\includegraphics[width=.5\linewidth]{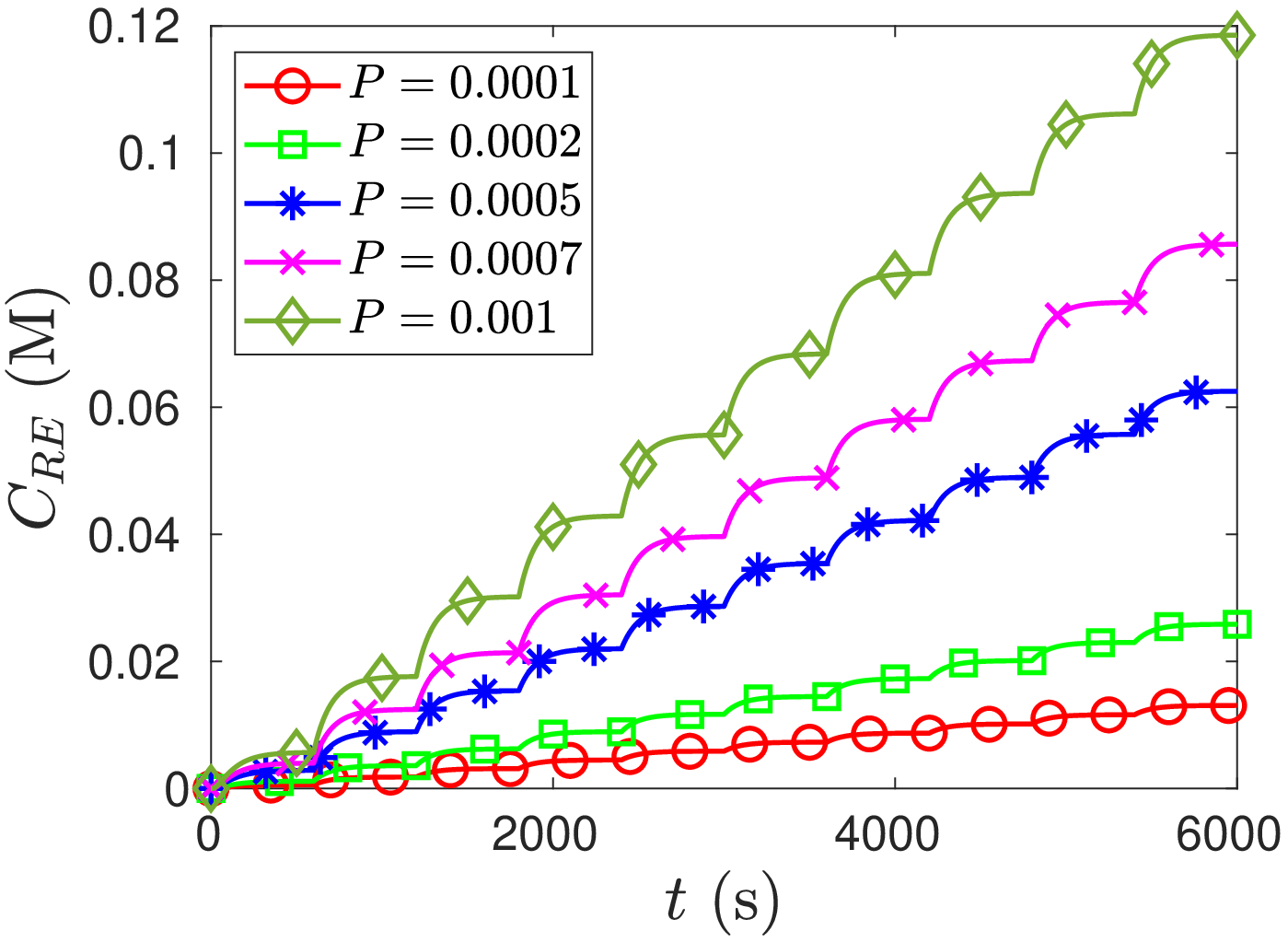}
		\\(a) $x=0.1$ mm & (b) $x=0.1$ mm \\	
		\includegraphics[width=.5\linewidth]{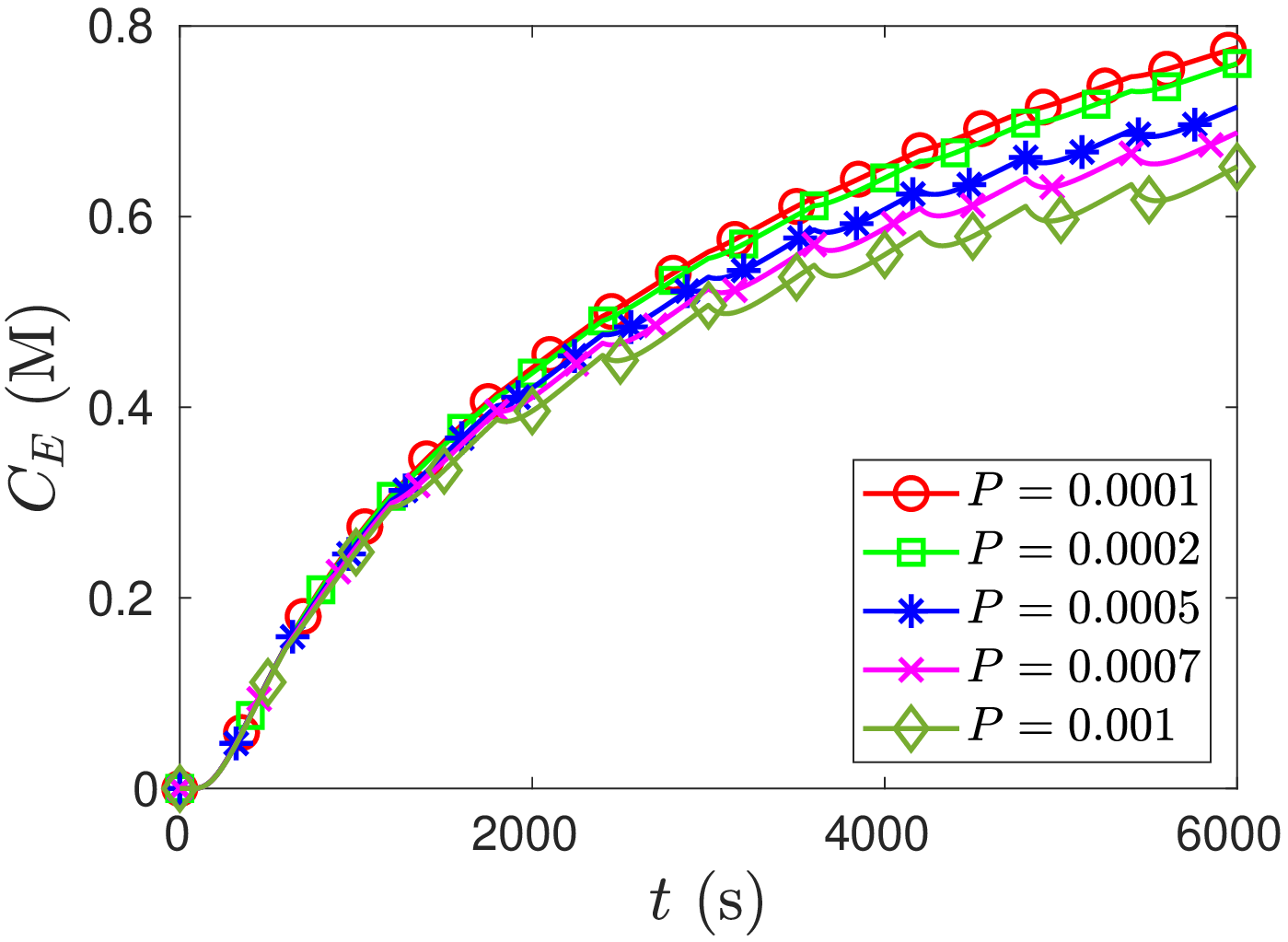}
		&
		\includegraphics[width=.5\linewidth]{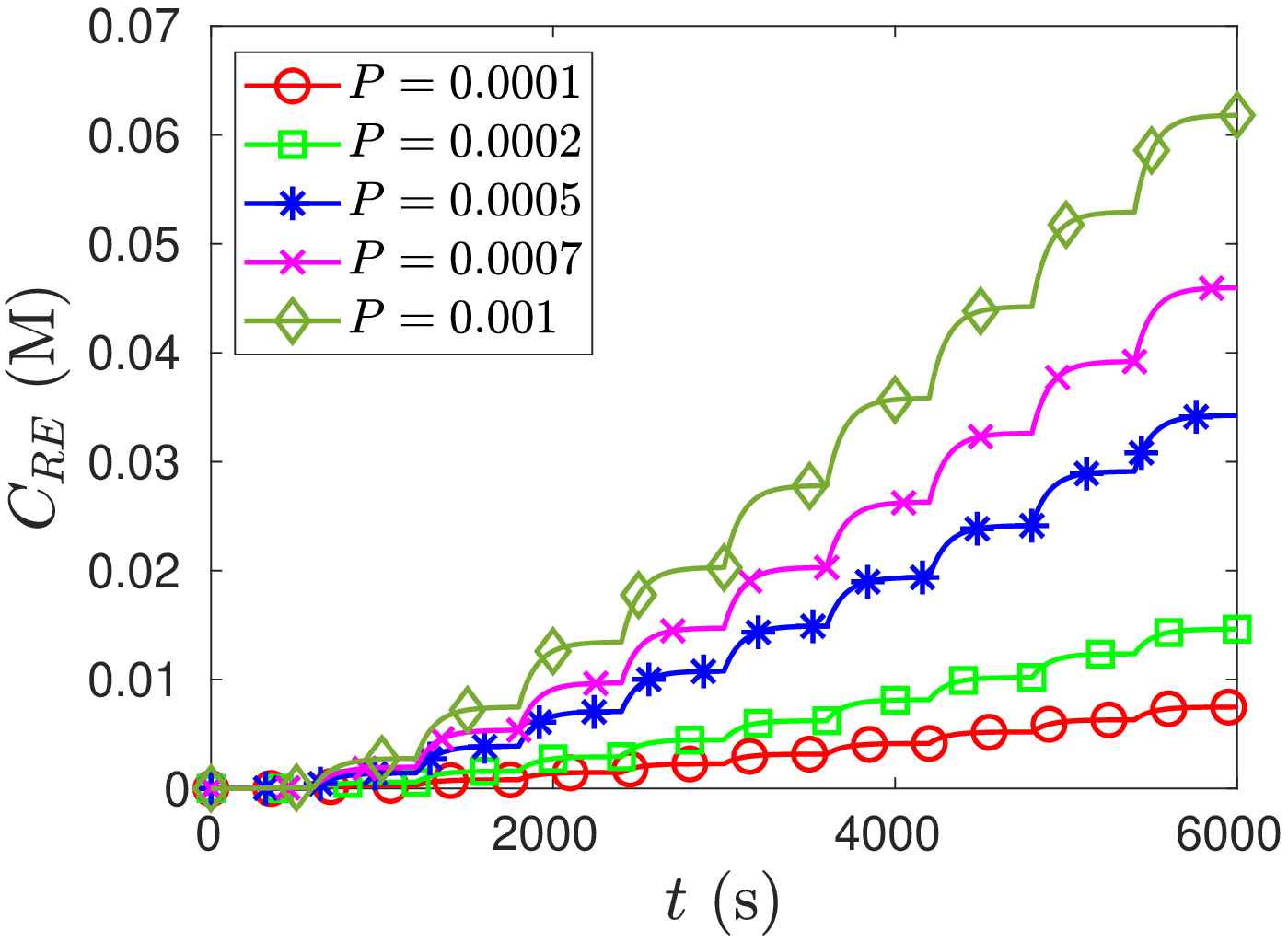}
		\\(c) $x=0.5$ mm & (d) $x=0.5$ mm \\
		\includegraphics[width=.5\linewidth]{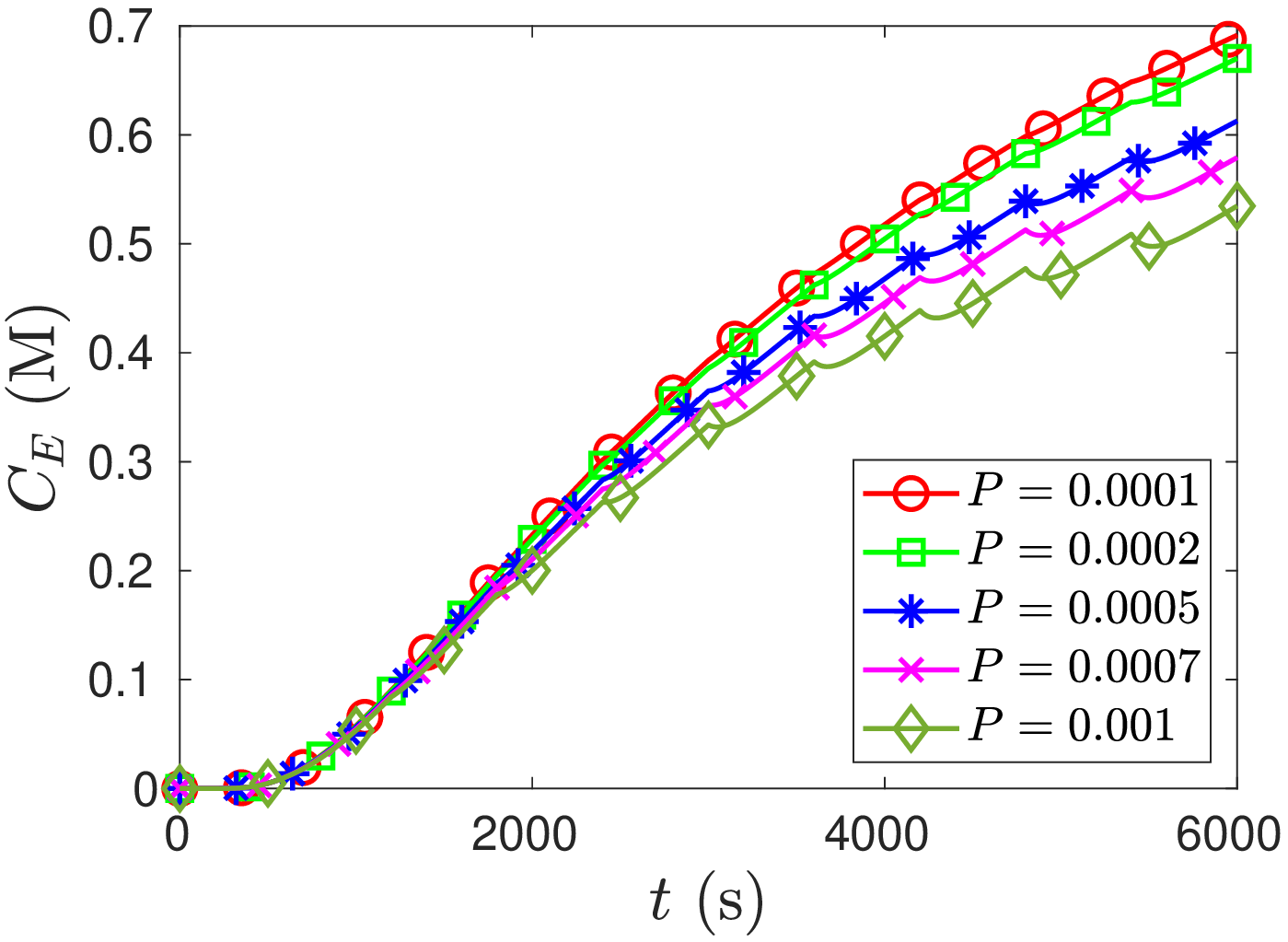}
		&
		\includegraphics[width=.5\linewidth]{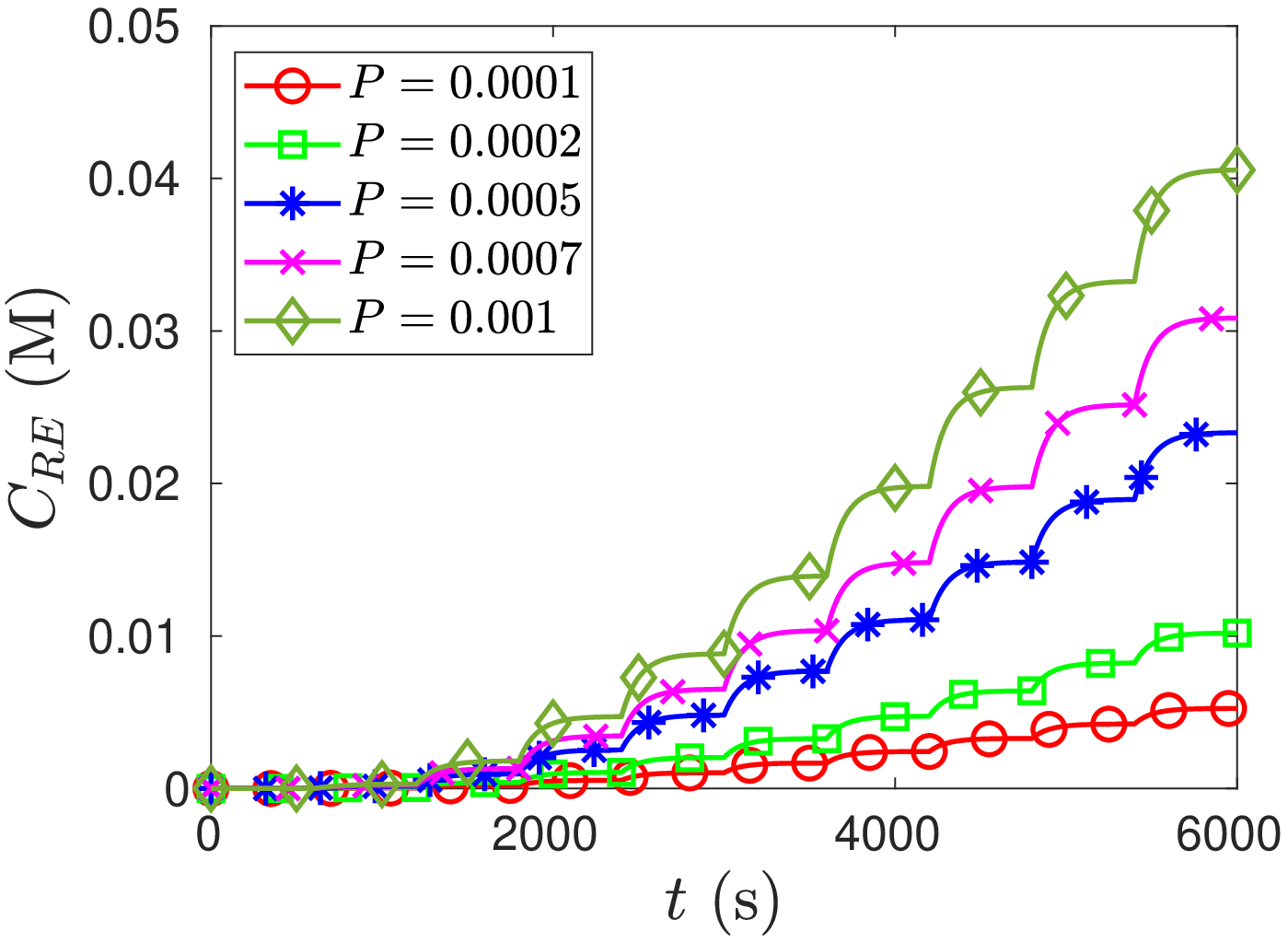}
		\\(e) $x=0.9$ mm& (f) $x=0.9$ mm\\		
	\end{tabular}
	\caption{Concentrations ($C_E$, $C_{RE}$) vs. time for $D=0.0001$ mm$^2$ s$^{-1}$ and different permeability (mm s$^{-1}$).}
	\label{2pchange}
\end{figure}

Fig. \ref{2pchange} illustrates the time-dependent drug transport at various locations throughout the tissue for different values of drug permeability. It is noticed from the  Fig. \ref{2pchange}a that at the initial stage, the extracellular concentrations increase rapidly as the point (0.1, 0.5) is close to the boundary where the drug is injected, so the drug diffusion to that point occurs quickly. Also, due to the same reason, the concentration profiles throughout the time are almost same for different permeability. But, in Figs. \ref{2pchange}c \& \ref{2pchange}f, it is observed that the concentrations for higher values of P are lesser than the concentrations for smaller P. The reason is that the drug uptake into the cells is more for higher P, which is shown in Figs.\ref{2pchange}d \& \ref{2pchange}f. Hence, the drug permeability enhances the drug uptake rate into the cells. One important observation is that after the application of ten pulses, the amount of drug entered into the cells is much lesser  compared to the amount of drug present in extracellular space. Slow mass transfer from extracellular to intracellular compartment is the possible reason.

\subsubsection{Effects of diffusion coefficient on drug concentrations}
\begin{figure}[h!]
	\centering
	\begin{tabular}{cc}
		\includegraphics[width=.5\linewidth]{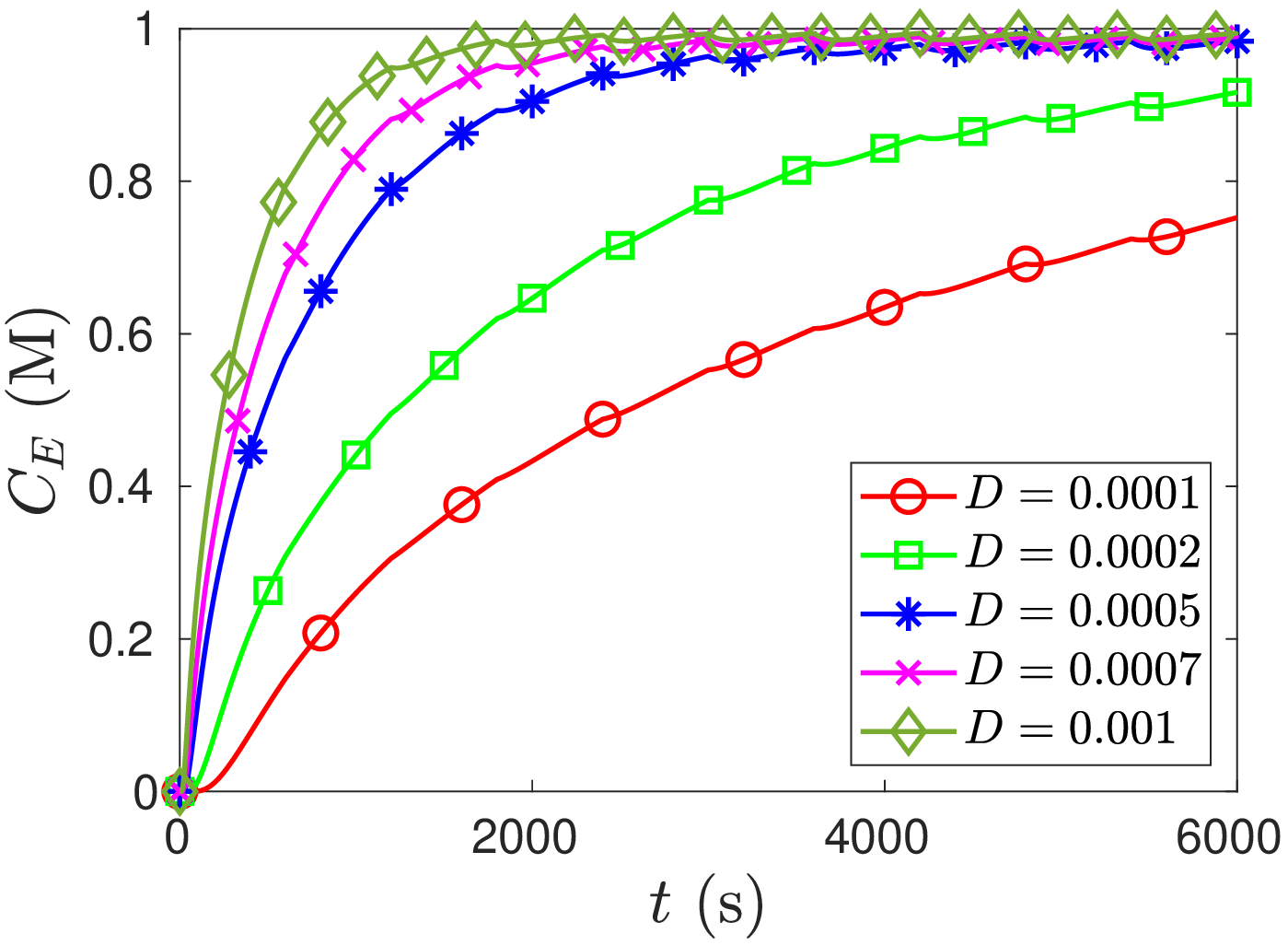}
		&
		\includegraphics[width=.5\linewidth]{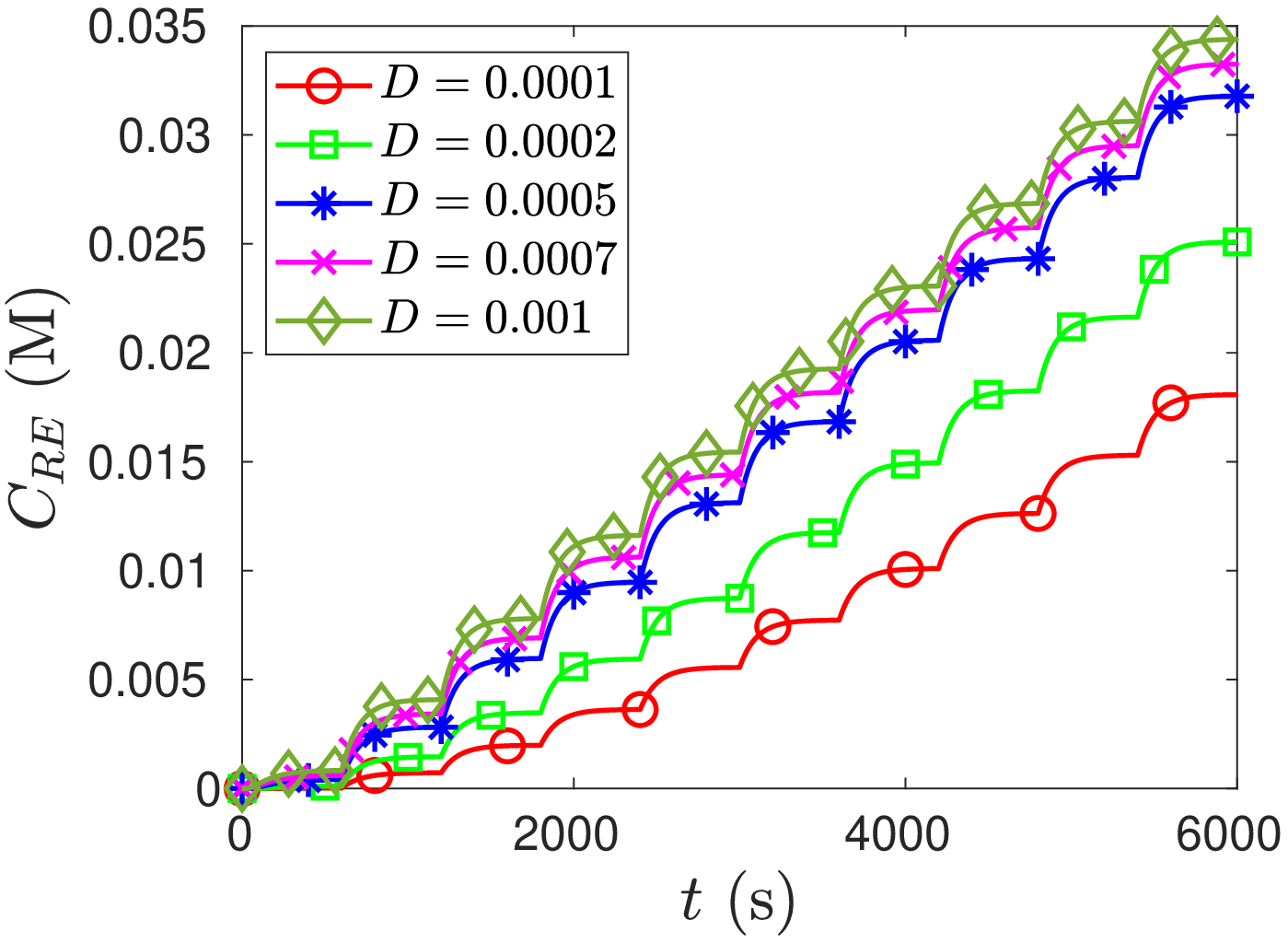}
		\\(a) & (b)
	\end{tabular}
	\caption{Concentrations ($C_E$, $C_{RE}$) vs. time for $P=0.0005$ mm s$^{-1}$  and different $D$ (mm$^2$ s$^{-1}$).}
	\label{2dchange}
\end{figure}
Fig. \ref{2dchange} represents the time-dependent drug concentration ($C_E$, $C_{RE}$) profiles for different diffusion coefficients (D) in the tissue at a point (0.5, 0.5). It is observed from the graphs that the drug concentrations ($C_E$, $C_{RE}$) increase with the increase in $D$ at all three locations. This is expected because, for a higher value of $D$, the drug particles move faster in the extracellular domain from the drug source and the drug uptake into cells is higher. As shown in Fig. \ref{2dchange}a, the concentrations of drug ($C_E$) have a sudden increase as  the drug diffusion  to the location (0.5, 0.5) occurs faster. After some time, the drug concentrations ($C_E$) increase steadily due to the drug uptake by the cells, which is shown in the Fig. \ref{2dchange}b.  As a result of faster drug diffusion in the ECS, intracellular drug concentrations ($C_{RE}$) also increase with the increase in $D$ that is noticed in Fig. \ref{2dchange}b. There are some jumps in the profiles of $C_{RE}$ because of repeated pulse application.

\subsubsection{Effects of pulse length and pulse number on drug concentrations}
\begin{figure}[h!]
	\centering
	\begin{tabular}{cc}
		\includegraphics[width=.5\linewidth]{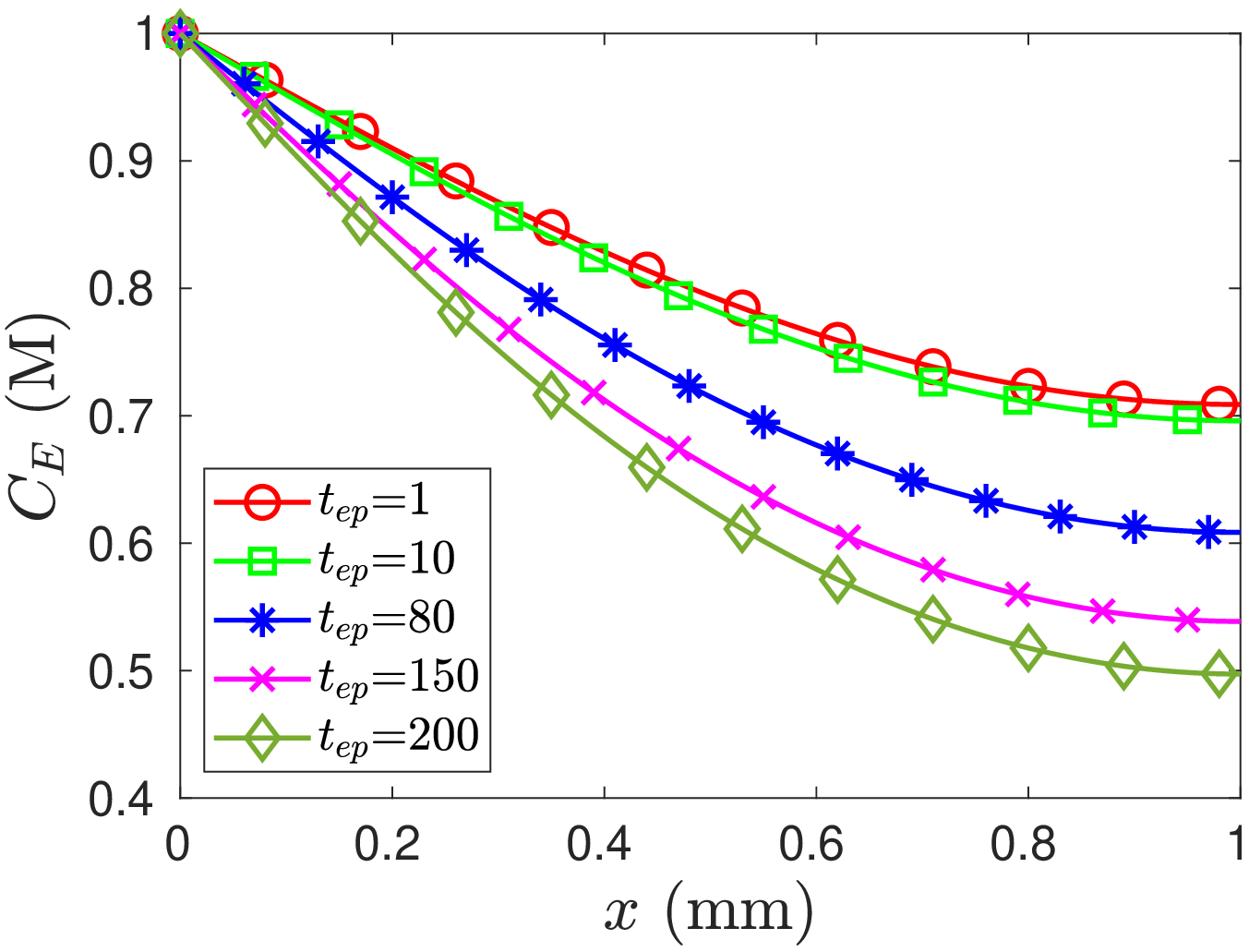}
		&
		\includegraphics[width=.5\linewidth]{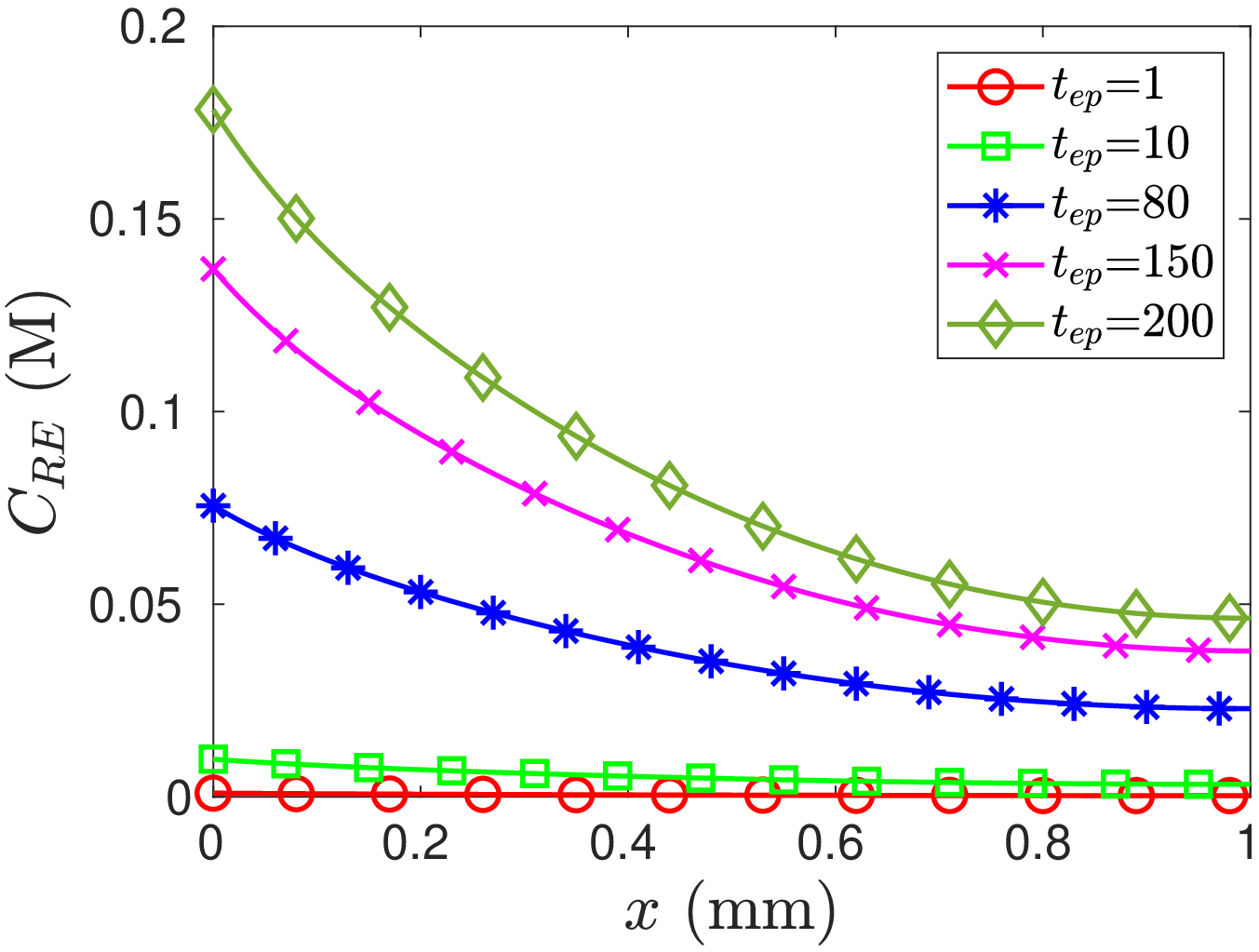}
		\\(a)  & (b)  \\	
	\end{tabular}
	\caption{The drug penetration along x-axis into the tissue after different pulse duration $t_{ep}$ ( ms). Concentrations ($C_E$, $C_{RE}$) are calculated for $D=0.0001$ mm s$^{-1}$ and $P=0.0005$ mm$^2$ s$^{-1}$. }
	\label{2tP_change}
\end{figure}
\begin{figure}[h!]
	\centering
	\begin{tabular}{cc}
		\includegraphics[width=.5\linewidth]{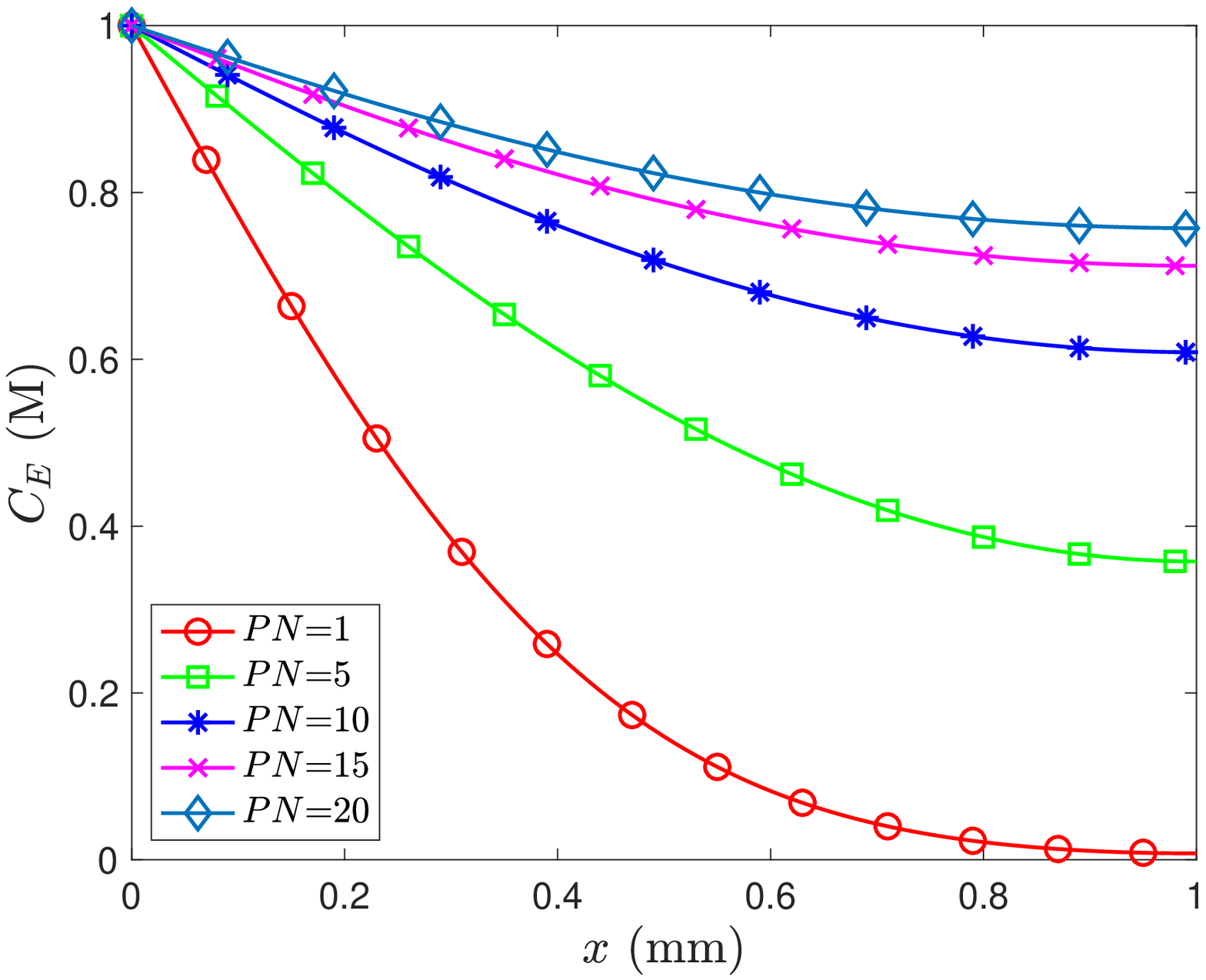}
		&
		\includegraphics[width=.5\linewidth]{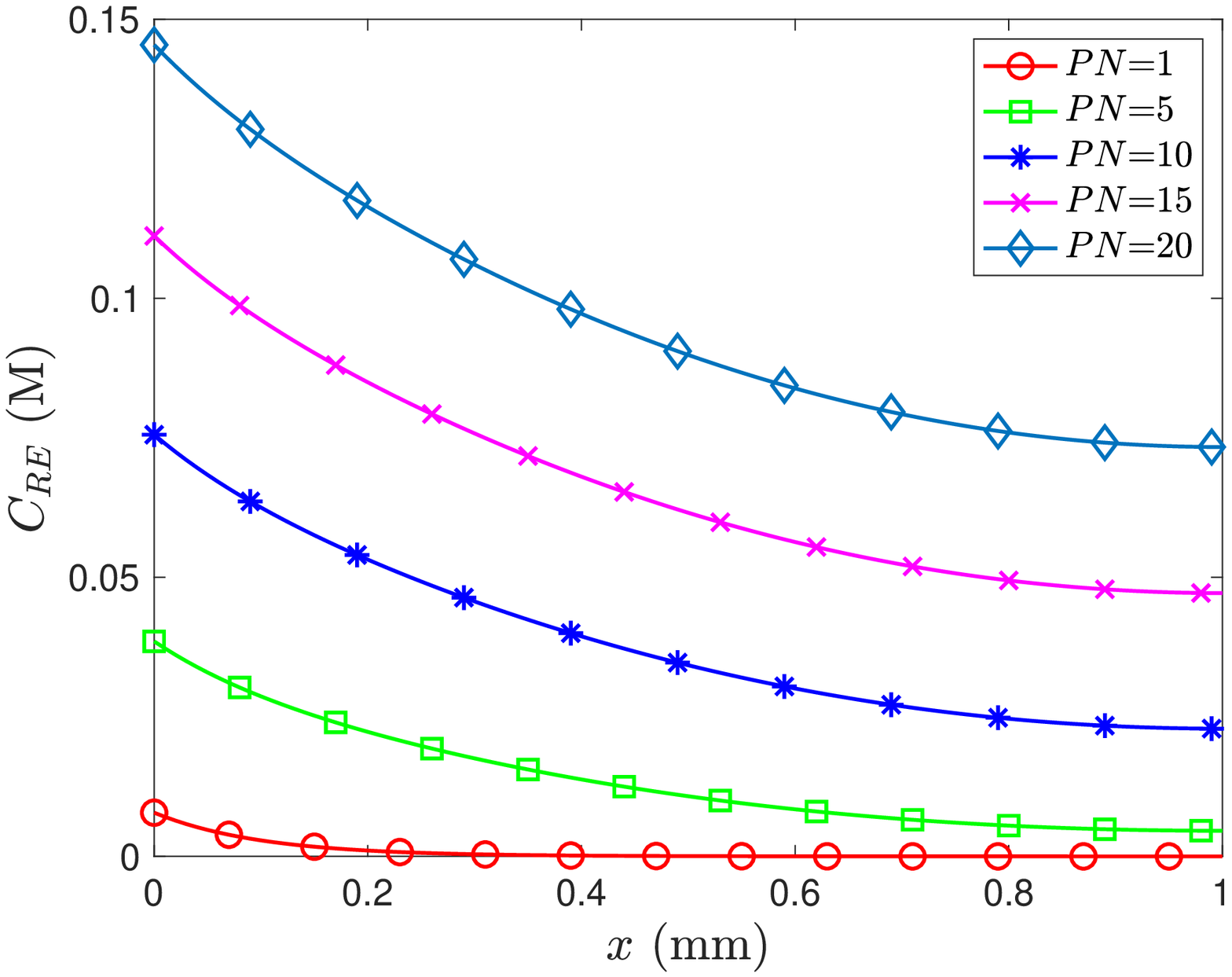}
		\\(a)  & (b) \\	
	\end{tabular}
	\caption{The drug penetration along x-axis into the tissue after different pulse number $PN$ of pulse length 80 ms. Concentrations ($C_E$, $C_{RE}$) are calculated for $D=0.0001$ mm s$^{-1}$ and $P=0.0005$ mm$^2$ s$^{-1}$.}
	\label{2PN_change}
\end{figure}
Fig. \ref{2tP_change} represents the drug concentration profiles in extracellular and intracellular compartments along the x-direction in the targeted tissue for different time duration of the applied pulses. Fig. \ref{2tP_change}a shows that the drug concentration in extracellular space decreases with the increase of pulse length. This is because the drug uptake into cells is higher for the application of longer pulses (see Fig. \ref{2tP_change}b). 
Also, it is noticed from the Fig. \ref{2tP_change}b that the drug uptake into the cells as well as the drug penetration throughout the tissue increases for the application of long duration pulses. The reason is that the pore density as well as mass transfer rate is increased with the application of long duration pulses, which is  already observed through the Fig. \ref{2mu}b.

Fig. \ref{2PN_change} shows the drug concentration profiles in extracellular and intracellular compartments along the x-direction in the tissue for different pulse numbers. It can be seen from both graphs that drug concentrations ($C_E$, $C_{RE}$) increase as $PN$ increases. $C_E$ increases as more drugs reach the ECS as a result of continuous diffusion for a long time. On the other hand, $C_{RE}$ increases as a result of the periodic reopening of membrane pores by the use of a pulse repeatedly and due to the availability of more drugs within the ECS.
Hence, a suitable number of a pulse that is applied periodically on the tissue can be determined as per the necessity of the amount of drug in the targeted cells.   

\subsubsection{Visualization of drug penetration by the contour plots}
\begin{figure}[h!]
	\centering
	\begin{tabular}{cc}
		\includegraphics[width=.5\linewidth]{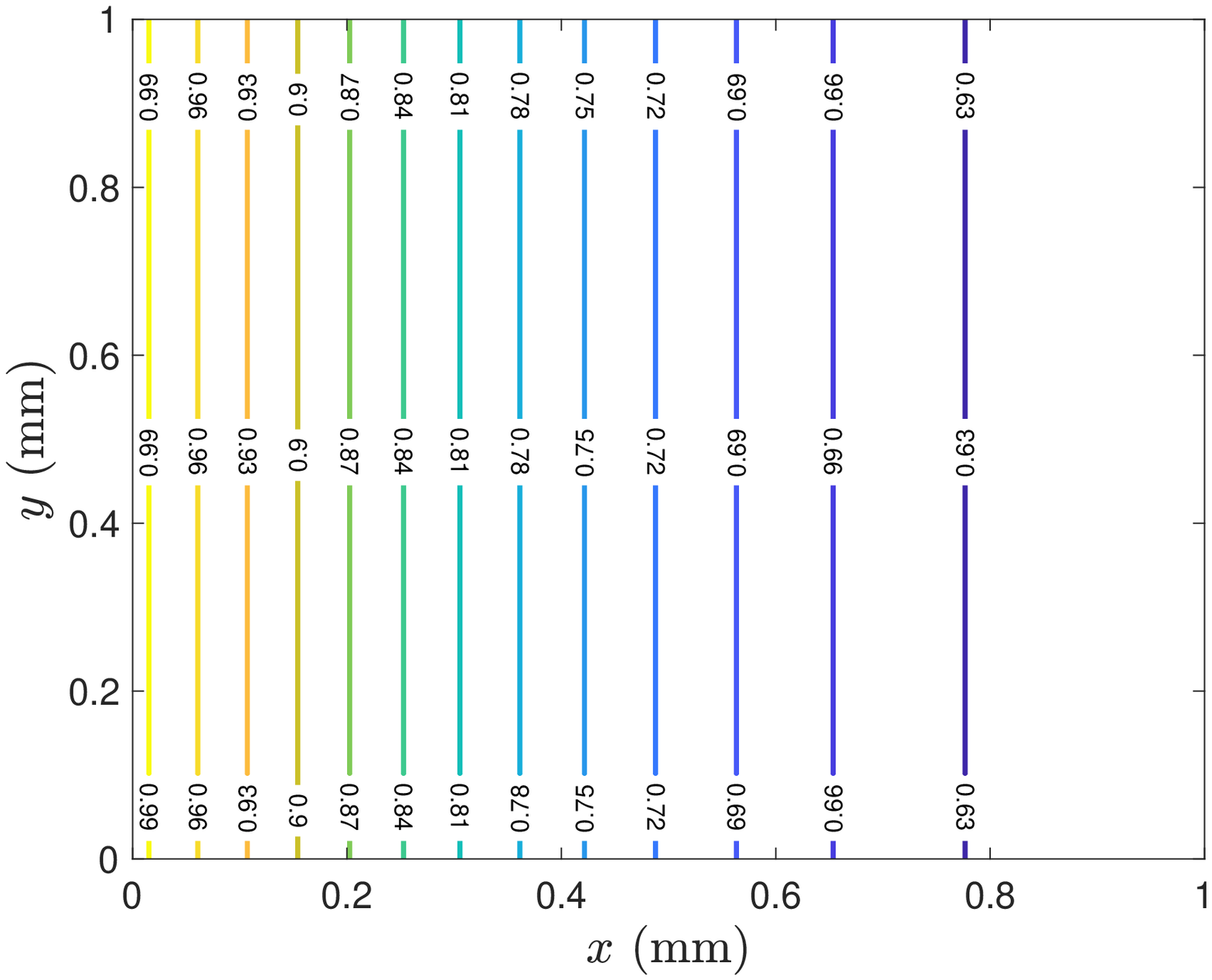}
		&
		\includegraphics[width=.5\linewidth]{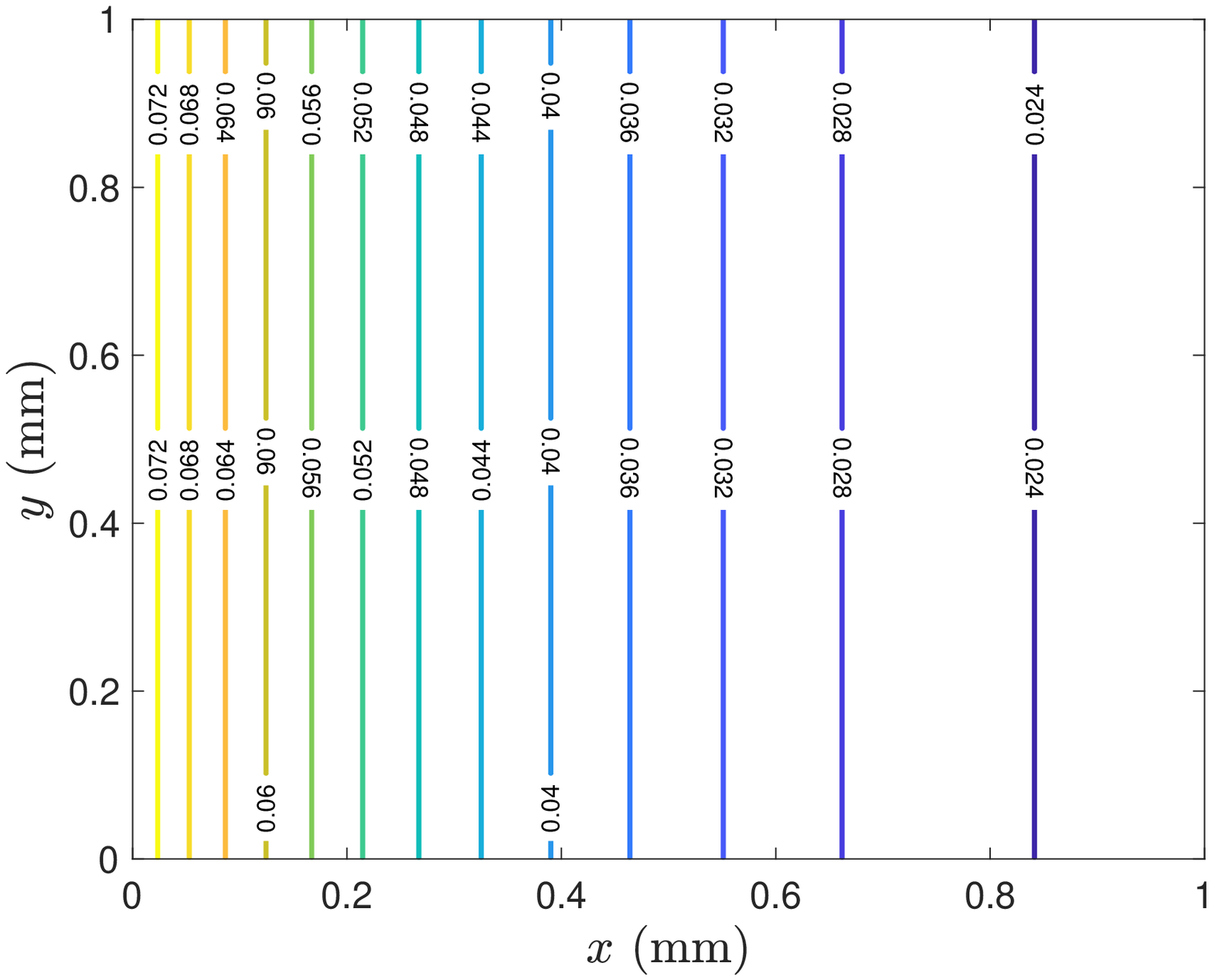}
		\\(a) $C_{E}$ & (b) $C_{RE}$ \\	
	\end{tabular}
	\caption{Contour plots of concentration distribution (iso-concentration lines) in the tissue. Here, $C_E$ and $C_{RE}$ are obtained for $P=0.0005$ mm s$^{-1}$, $D=0.0001$ mm$^2$ s$^{-1}$, $PN=10$ and $t_{ep}=80$ ms.}
	\label{2contour}
\end{figure}
The drug penetration throughout the targeted tissue after electroporation with ten pulses is shown by the contour plots (Fig. \ref{2contour}). It is seen that after performing the electroporation of ten pulses with time gap of 600  s, a large amount of drug enters into the cells near the left boundary, whereas a lesser amount of drug reaches the cells near the right boundary. It takes different times for the drug to enter the cell at different locations in the tissue. These times are calculated numerically and given in the Table \ref{2tab2}.
Therefore, in this model's approach, a dose of drug with a concentration of more than 0.025 M can be administered into cells of a targeted tissue. This amount of drug could be a desired requirement for the treatment of the tissue. However, if more drug is required to treat a deceased tissue, it will be possible to increase the amount of drug into the cells by applying more pulses.

\begin{table}[h!]
	\caption{Time taken to initiate drug into the cell at various locations are shown.}
	\vspace{0.2cm}
	\label{2tab2}
	\centering
	\begin{tabular}{|c|c|c|c|}
		\hline 
		& &  \\ [-1em]
		Location & Time (s)  & $C_{RE}$ (M)  \\ [0.2em]
		\hline 
		& &  \\ [-1em]
		(0.1, 0.5)& 2.0$\times 10^{-4}$ & 1.614850$\times 10^{-56}$\\ [0.2em]
		\hline 
		& &  \\ [-1em]
		(0.2, 0.5)& 4.0$\times 10^{-4}$ & 1.653603$\times 10^{-103}$\\ [0.2em]
		\hline 
		& &  \\ [-1em]
		(0.3, 0.5)& 6.0$\times 10^{-4}$ & 1.693286$\times 10^{-150}$\\ [0.2em]
		\hline 
		& &  \\ [-1em]
		(0.4, 0.5)& 8.0$\times 10^{-4}$ & 1.733921$\times 10^{-197}$\\ [0.2em]
		\hline 
		& &  \\ [-1em]
		(0.5, 0.5)& 1.0$\times 10^{-3}$ & 1.775532$\times 10^{-244}$\\ [0.2em]
		\hline
		& &  \\ [-1em]
		(0.6, 0.5)& 1.2$\times 10^{-3}$ & 1.818141$\times 10^{-291}$\\ [0.2em]
		\hline
		& &  \\ [-1em]
		(0.7, 0.5)& 1.4$\times 10^{-3}$ & 1.861773$\times 10^{-338}$\\ [0.2em]
		\hline
		& &  \\ [-1em]
		(0.8, 0.5)& 1.6$\times 10^{-3}$ & 1.906451$\times 10^{-385}$\\ [0.2em]
		\hline
		& &  \\ [-1em]
		(0.9, 0.5)& 1.8$\times 10^{-3}$ & 1.952202$\times 10^{-432}$\\ [0.2em]
		\hline
	\end{tabular}
\end{table}

\subsection{Thermal effects in the tissue}
\begin{figure}[h!]
	\centering
	\begin{tabular}{cc}
		\includegraphics[width=.5\linewidth]{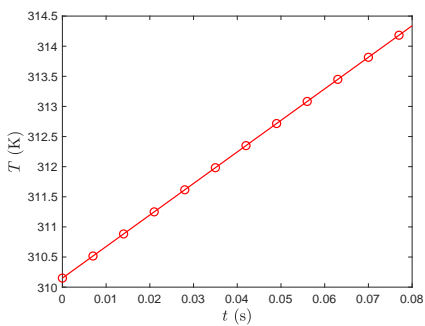}
		&
		\includegraphics[width=.5\linewidth]{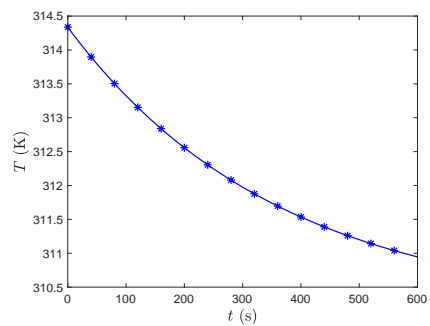}
		\\(a) During electroporation  & (b) After electroporation \\	
	\end{tabular}
	\caption{The change of temperature $T$ with time $t$ at (0.5, 0.5).}
	\label{2temp_time}
\end{figure}
Fig. \ref{2temp_time}a shows the change of temperature at the middle (0.5, 0.5) of the tissue during electroporation due to the application of one pulse of 80 ms. It is clearly observed that the temperature increases with the pulse time from the initial body (tissue) temperature 310.15 K due to Joule heating. However, Fig. \ref{2temp_time}b describes the temporal temperature variations at (0.5, 0.5)  after electroporation. From the figure, it is seen that the temperature gradually decreases with time from its peak value ($T_{ep}=314.335$ K), and after a certain time around 600  s, it reaches the temperature ($\approx 311$ K) that is close to the body temperature $T_{b}$.  This happens due to the continuous heat transfer from  inside to the outside of the tissue.

\begin{figure}[h!]
	\centering
	\includegraphics[width=.5\linewidth]{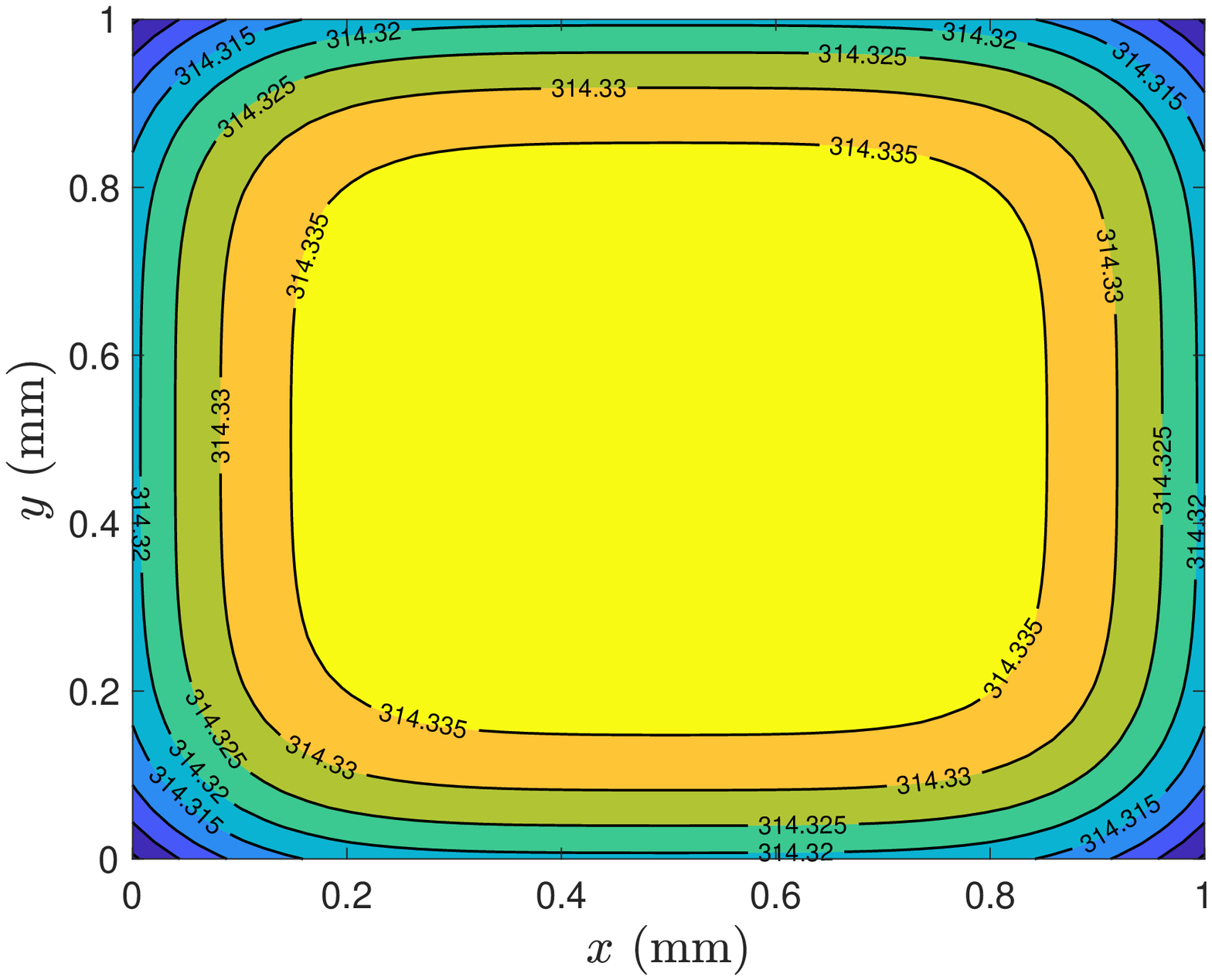}
	\caption{Temperature distribution over the tissue region after electroporation with one pulse.}
	\label{2temp_surface} 
\end{figure}
Fig. \ref{2temp_surface} illustrates the temperature contours throughout the tissue region after the electroporation over a time duration of 80 ms. It is noticed that the maximum temperature is attained at the center, and the temperature becomes minimum at the corner points of the tissue region. The maximum value of the increased temperature is about 314.335 K, which is less than the possible temperature of  cell damage (315.15 K $=42^\circ C$)  as reported in \cite{Garcia2014}. So, there is no chance of cell damage or cell death in the tissue due to application of the pulses. 
The main reason for slight increase in temperature due to electroporation is the application of low voltage pulse. One important observation is that in order to  avoid any cell damage, the maximum pulse duration should be of 80 ms for the application of the electric field 28 V mm$^{-1}$.
Hence, the electroporation configuration and model arrangement are options for successful drug delivery into the targeted tissue with reversible electroporation.

\section{Conclusions} 
From this study,  the following observations are made:
\begin{itemize}
	\item The proposed mathematical model for drug delivery deals with an advanced electroporation technique. Low voltage multiple pulses are applied to introduce the drug into the cells of a diseased tissue with no cell death.
	\item Since, the initial drug is injected from one side (left boundary) of the tissue, it takes a long time to spread the drug into the targeted tissue.  
	
	\item The drug uptake into the cells increases with the pulse duration and the pulse number.
	
	\item It is observed that during electroporation, the tissue temperature increases from its initial body temperature (310.15 K) due to Joule heating. During post electroporation period, the temperature decreases with time from its peak value and it tends to reach the normal body temperature.
	
	\item In this model, no cell damage  takes place due to thermal effects as increased temperature is not that significant.
 
	\item The pulse length (80 ms) and the rest time between the pulses (600  s) are optimized to deliver a desired amount of drugs into the cells with no thermal damage. This is a novel outcome of this model.
\end{itemize}
The direct experimental validation of the model is not possible due to lack of necessary data. However, the model is compared qualitatively with experimental and theoretical studies. The model could be useful in various clinical experiments and pharmaceutical sectors.

%
%
%
%

\bmhead{Acknowledgments}

The first author (Nilay Mondal) acknowledges the support received from the Indian Institute of Technology Guwahati, Assam, India.

 \section*{Conflict of interest}

On behalf of all authors, the corresponding author states that there is no conflict of interest.

\bibliography{sn-article.bbl}



\end{document}